\documentclass[12pt,thmsa]{article}
\usepackage{amsmath}
\usepackage{graphicx}
\usepackage{amsfonts}
\usepackage{amssymb}
\usepackage{mathrsfs}
\newtheorem{theorem}{Theorem}[section]
\newtheorem{lemma}[theorem]{Lemma}

\newtheorem{remark}{Remark}

\newtheorem{Theo}{Theorem}[section]

\begin{document}

\begin{center}
{\Large \textbf{Smoothed functional average variance estimation    for dimension reduction}}

\bigskip

M\`etolidji Moquilas Raymond AFFOSSOGBE\textsuperscript{a} , Guy Martial  NKIET\textsuperscript{b}  and Carlos OGOUYANDJOU\textsuperscript{a}

\bigskip

\textsuperscript{a}Institut de Math\'ematiques et de Sciences Physiques, Porto Novo, B\'enin.
\textsuperscript{b}Universit\'{e} des Sciences et Techniques de Masuku,  Franceville, Gabon.

\bigskip

E-mail adresses : metolodji.affossogbe@imsp-uac.org,  guymartial.nkiet@mathsinfo.univ-masuku.com,  ogouyandjou@imsp-uac.org.

\bigskip
\end{center}

\noindent\textbf{Abstract.}We propose an estimation method that we call functional average variance estimation (FAVE), for estimating the EDR space in functional semiparametric regression model,  based on kernel estimates of density and regression. Consistency results are then established for the estimator of the interest operator, and for the directions of EDR space. A simulation study that shows that the proposed approach performs as well as traditional ones is presented.

\bigskip


\section{Introduction}
\label{Intro}
In recent years, much  attention has been given to functional statistics, which can be described as the set of statistical methods for processing data having  the form of curves considered as observations of functions belonging to given functional spaces.   Among the references  in this field, there  are the books  by Ramsay and Silverman\cite{ramsay}  for the applied aspects, Bosq \cite{bosq} for the theoretical aspects,  Ferraty and Vieu \cite{ferraty}, and Horv\'ath and Kokoszka\cite{horvath} for recent developments. Many works in this field deal with problems   that appear  in the general framework of functional regression models which  are usually used to find the best link between a real random variable $Y$  and a random curve $X$ whose values belong to $\mathcal{ H} = L ^ 2 ([0,1])$, the set of square integrable functions from $ [0,1]$ to  $ \mathbb {R} $. An abundant literature has examined cases of parametric functional regression models (e.g., \cite{cardot}, \cite{ramsay}, \cite{hall}, \cite{yao})  described by the relation $Y = f_{\theta} (X, \varepsilon) $, where $ f_{\theta} $  belongs to a well-known family of functions parameterized by the unknow parameter $\theta$  which is to be estimated, and $\varepsilon$ is an  error term. 
In contrast to this, some works deal  with a nonparametric model  $Y=f(X)+\varepsilon $ where $f$ is an unknown and arbitrary function to estimate, and have introduced nonparametric estimation approaches, such as methods  based on kernel estimators (\cite{ferraty},\cite{ferraty2}). Alternatively, between these two different approaches, a semiparametric regression model 
\begin{align}\label{reg}
Y=f(<\beta_1,X>_\mathcal{H}, <\beta_2,X>_\mathcal{H},\cdots , <\beta_K,X>_\mathcal{H},\varepsilon)
\end{align}
was considered (\cite{dauxois},\cite{ferre},\cite{ferre2},\cite{lian}). In the model (\ref{reg}), $<\cdot,\cdot>_\mathcal{H}$ denotes the inner product of $\mathcal{H}$ defined for all $g_1$ and $g_2$ belonging to  $ \mathcal{H} $ by $<g_1, g_2>_\mathcal{H} = \int_{0 }^{1} g_{1}(t)g_{2}(t)dt$, and $\beta_1,\cdots,\beta_K$ are elements of $\mathcal{H}$ to be estimated. This model just is an extension in the functional case of the model introduced by Li\cite{li} in the multivariate context and which has been intensively studied since then. It expresses the fact that the information in $X$  about  $Y$  depends only on
the projection of $X$ onto the subspace  spanned by
$\left\{\beta_{1},\cdots,\beta_{K}\right\}  $, called  effective dimension-reduction
(EDR) space. Li\cite{li} showed that the problem of estimating the EDR  space  comes down, under a fairly general condition, to the
spectral analysis of an  operator depending on  the covariance
operator of the conditional expectation $\mathbb{E}\left(  \left.  X\right|
Y\right)   $ of $X$ given $Y$. Then,  he proposed an estimation method, called sliced inverse regression (SIR),  based on an estimate of an approximation of this covariance operator  obtained by slicing the range of $Y$. Alternatively, Cook\cite{cook} proposed another method, called sliced average variance estimation (SAVE), for estimating the EDR by using  an estimate of an approximation of an operator depending on the conditional covariance operator $Var(X\vert Y)$ of $X$ given $Y$. SIR and SAVE are the most popular methods for dimension reduction in the multivariate context, and smoothed estimaton methods, based on kernel estimates,  have been proposed  for them respectively by Zhu and Fang\cite{zhu} and Zhu and Zhu\cite{zhu2}. In the functional context, SIR has been extended to functional SIR (FSIR)  by Ferr\'e   and Yao \cite{ferre} who also  proposed later a smoothed estimation procedure based on kernel estimates, so defining  smoothed functional inverse regression (FIR). On the other hand, more recently, Lian and Li\cite{lian} extended SAVE to functional SAVE (FSAVE). To the best of our knowledge, a smoothed estimation of  SAVE have not been proposed yet in the context of functional data. Taking all this  into consideration, we introduce  in this paper a  kernel functional average variance estimation  (FAVE) method for estimating the EDR space related to model (\ref{reg}). The rest of the paper is organized as follows.  In Section 2, we recall some basic facts about FAVE in the functional context, and we specify the interest operator to estimate. In Section 3, an estimator based on kernel estimates is proposed for this estimating this operator. Section 4 is devoted to an asymptotic study of the introduced estimator. A simulation study that permits to evaluate the performance of our proposal is presented in Section 5. The proofs of theorems are postponed in Section 6.

\section{ Functional Sliced Average Variance}\label{sec2}
Let us  consider the random variables $Y$ and $X$ involved in the model (\ref{reg});  we assume, without loss of generality, that  $\mathbb{E}(X)=0$, and that  $\mathbb{E}(\Vert X\Vert_\mathcal{H}^2)<+\infty$. Then, the   covariance operator of $X$ is defined by $\Gamma=\mathbb{E}(X\otimes X),$ where for any $x,y\in\mathcal{ H}$, $x\otimes y$ denotes the linear operator from $\mathcal{ H}$ to itself such that $ (x\otimes y)(h)=<x,h>_\mathcal{H}y$ for any $h\in\mathcal{H}$. Throughout the paper, $\Gamma$ will be assumed to be non-singular and positive definite. Letting $\mathcal{B}=(<\beta_1,X>_\mathcal{H}, <\beta_2,X>_\mathcal{H},\cdots , <\beta_K,X>_\mathcal{H})$ and denoting by $Var(X\vert\mathcal{B})$ the conditional covariance operator of $X$ given $\mathcal{B}$,  we consider the following assumptions:

\bigskip

\noindent $(\mathscr{A}_1)$ :  for all $b\in\mathcal{H}$, one has $\mathbb{E}(<b,X>_\mathcal{H}|\mathcal{B})=\sum\limits_{k=1}^{K}c_k<\beta_k,X>_\mathcal{H}$, where $c_1,\cdots,c_K$ are real numbers;

\bigskip

\noindent  $(\mathscr{A}_2)$: $Var(X|\mathcal{B})$ is nonrandom.

\bigskip
 
 \noindent Lian and Li\cite{lian} showed that   under the  assumptions $(\mathscr{A}_1)$ and   $(\mathscr{A}_2)$, one has  the inclusion 
\begin{equation}\label{inclus}
R\left(\Gamma-Var(X|Y)\right)\subset\Gamma \Im , 
\end{equation}
where $R(A)$ denotes the range of the operator $A$,  $Var(X\vert Y)$ denotes the conditional covariance operator of $X$ given $Y$, and $\Im$ is the EDR space, that is the space spanned by $\beta_1,\cdots,\beta_K$. 
Therefore, $R(\Gamma_I)\subset\Im$, where
\[
\Gamma_I:=\Gamma^{-1}\mathbb{E}\bigg(\Gamma-2Var(X|Y)+Var(X|Y)\Gamma^{-1}Var(X|Y)\bigg).
\] 
An important consequence is that $\Gamma_I$ is degenerate in any direction orthonormal to the $\beta_k$'s ($k=1,2,\cdots,K$). Then  $\Gamma_I$ is a finite rank operator whose range is contained into the EDR space. This  space can, therefore,  be approached by the subspace spanned by the eigenvectors of $\Gamma_I$ associated with the $K$ largest non-null eigenvalues of $\Gamma_I$ in the same way as in the multivariate case. In the following we suppose that rank($\Gamma_I$)=K. We see, therefore, that the eigenvectors associated with the $K$ largest eigenvalues of $\Gamma_I$ form an base to EDR space, which make the EDR space identifiable. So $\Gamma_I$ is the interest operator of the FAVE method.   Since the domain of $\Gamma^{-1}$ is not the whole $\mathcal{H}$, $\Gamma_I$ may not be well-defined. Conditions under which this operator is well defined are established in \cite{lian} and recalled below.

 \bigskip
 
 \noindent Let $$X=\sum_{j=1}^{\infty}\xi_j\phi_j,$$ be the well-known Karhunen-Loève expansion of $X$, where $\mathbb{E}[\xi_{j}^{2}]=\alpha_j$ are the eigenvalues and $\phi_j$ are the eigenfunctions. As usual in the functional data literature (e.g.,\cite{lian},\cite{ferre2}), we assume that  $\alpha_1>\alpha_2>\cdots>0$. We now introduce the assumptions:

\bigskip

			\noindent $(\mathscr{A}_3)$: $\mathbb{E}\left( \Vert X\Vert_{\mathcal{H}}^{4}\right)<+\infty $ ;

\bigskip

			\noindent $(\mathscr{A}_4)$:  $\mathbb{E}\left[ \left(\sum\limits_{j=1}^{\infty}\alpha_{j}^{-2}\sum\limits_{i=1}^{\infty}Cov^{2}(\xi_i,\xi_j|Y)\right) ^2\right]<+\infty.$ 

\bigskip

\noindent It is known that if  $(\mathscr{A}_3)$ and $(\mathscr{A}_3)$ hold, then $\Gamma_I$ is well-defined  (see Proposition 1 in  \cite{lian}).

\section{ Kernel estimator of the interest operator }\label{sec2}

For performing the FAVE method $\Gamma_I$,  has to be estimated. Lian and Li\cite{lian} introduced an estimator obtained by slicing the range of $Y$. In this section, we propose another estimator of this operator based on kernel estimates  of  density and   regression. Since  $\Gamma=\mathbb{E}[Var(X|Y)]+Var[\mathbb{E}(X|Y)]$, we have  
\begin{equation}\label{gammai}
\Gamma_I=\Gamma^{-1}\bigg(2\Gamma_e+\Psi-\Gamma\bigg)
\end{equation}
where $\Gamma_e= Var\left[ \mathbb{E}(X|Y)\right]$ is the covariance operator of the conditional expectation  $\mathbb{E}(X|Y)$ and  $\Psi= \mathbb{E}\left( Var(X|Y)\Gamma^{-1}Var(X|Y)\right)$. Ferr\'e and Yao \cite{ferre} introduced a  kernel estimator of $\Gamma_e$ and showed  its consistency. Here, we will use this estimator, and also a kernel  estimator  of  $\Psi$ together with the empirical counterpart of $\Gamma$ in order to define an estimator of $\Gamma_I$. Letting $f$ be the density of $Y$ and putting
\[
m(y)=\mathbb{E}(\textbf{1}_{\{Y=y\}}\,X),\,\,\,M(y)=\mathbb{E}(\textbf{1}_{\{Y=y\}}\,X\otimes X),
\]
\[
r(Y)=\mathbb{E}(X|Y)=\dfrac{m(Y)}{f(Y)}\,\,\,\textrm{ and }\,\,\,R(Y)=\mathbb{E}(X\otimes X|Y)=\dfrac{M(Y)}{f(Y)},
\]
we have  $\Psi=\mathbb{E}\left( C(Y)\Gamma^{-1}C(Y)\right)$ where  $C(Y)=Var(X|Y)=R(Y)-r(Y)\otimes r(Y)$. As it was done in \cite{zhu}, in order to
avoid the effect of the small values in the denominator, we consider     $f_{e_n} = \max(f, e_n)$ instead of $f$, where $(e_n)_{n\in \mathbb{N}^{*}}$ is  a sequence of real numbers which tends to zero as $n\rightarrow +\infty$. Then, we consider
\[
r_{e_n}(Y)=\dfrac{m(Y)}{f_{e_n}(Y)},\,\,\, R_{e_n}(Y)=\dfrac{M(Y)}{f_{e_n}(Y)}\,\,\,\textrm{ and }\,\,\,C_{e_n}(Y)=R_{e_n}(Y)-r_{e_n}(Y)\otimes r_{e_n}(Y)
\]
instead of $r(Y)$, $R(Y)$ and $C(Y)$. The definition of $\Gamma_I$ given in (\ref{gammai}) requires to use the inverse of $\Gamma$. But since 
	$\Gamma$ is an Hilbert-Schmidt operator, even though its  inverse  exists  it is not generally bounded. To avoid this difficulty, we consider instead  the finite-rank operator $\Gamma_{D}=\Pi_{D}\Gamma\Pi_{D}$, where $D\in\mathbb{N}^\ast$ and  $\Pi_{D}$ is the  projector onto the subspace $S_{D}$  spanned by the system $\{\phi_1,\cdots,\phi_D\}$ consisting of the $D$  first elements  of an orthonormal  basis of $\mathcal{H}$. This basis can, for example,  be obtained either from principal component analysis (PCA)  of $X$ or by using $B$-splines basis.  This operator has a   bounded (pseudo-)inverse defined by $\Gamma_{D}^{-1}=\Pi_{D}\Gamma^{-1}\Pi_{D}$. 

\bigskip

\noindent Let  $\{(X_i, Y_i)\}_{1\leq i\leq n} $ be an i.i.d. sample of   $(X,Y)$; the empirical counterpart of $\Gamma$ is given by $\widehat{\Gamma}_{n}={n}^{-1}\sum\limits_{i=1}^{n}X_i\otimes X_i$. Considering the estimate $ \widehat\Pi_{D}$ of $\Pi_D$ defined as the projector onto an estimate $\widehat{S}_D$ of  $S_D$, we estimate   $\Gamma_{D}$ by $  \widehat{\Gamma}_{D}= \widehat\Pi_{D}\widehat{\Gamma}_{n} \widehat \Pi_{D}$. If we use PCA (resp. $B$-splines basis)  then $\widehat{S}_D$ consists of the eigenvectors associated with the $D$ largest eigenvalues of $\widehat{\Gamma}_n$ (resp. $\widehat{S}_D=S_D$).
For a given kernel function  $K\,:\,\mathbb{R}\rightarrow\mathbb{R}_+$ and a given real $h>0$,   we consider the estimates
\[
  \widehat{f}(y)= \dfrac{1}{nh}\sum\limits_{i=1}^{n}K\left(\dfrac{Y_i-y}{h} \right),\,\,\,\widehat{m}(y)=\dfrac{1}{nh}\sum\limits_{i=1}^{n}K\left(\dfrac{Y_i-y}{h} \right)\, X_i
\]
and
\[
\widehat{M}(y)=\dfrac{1}{nh}\sum\limits_{i=1}^{n}K\left(\dfrac{Y_i-y}{h} \right)\,\,X_i\otimes X_i
\]
of  $f$, $m$ and $M$ respectively. Then, putting
\[
\widehat{f}_{e_n}(y)=\text{max}\{e_n,  \widehat{f}(y)\},\,\,\,\widehat{r}_{e_n}(y)=\dfrac{ \widehat{m}(Y)}{ \widehat{f}_{e_n}(y)},\,\,\,
 \widehat{R}_{e_n}(y)=\dfrac{ \widehat{M}(y)}{ \widehat{f}_{e_n}(y)}
\]
and 
\[
\widehat{C}_{e_n}(y)= \widehat{R}_{e_n}(y)- \widehat{r}_{e_n}(y)\otimes  \widehat{r}_{e_n}(y)
\]
we consider
\[
\widehat{\Gamma}_{e,n}=\dfrac{1}{n}\sum_{i=1}^{n} \widehat{r}_{e_n}(Y_i)\otimes  \widehat{r}_{e_n}(Y_i),\,\,\,
\widehat{\Psi}_{e_n,D}=\dfrac{1}{n}\sum\limits_{i=1}^{n} \widehat{C}_{e_n}(Y_i) \widehat{\Gamma}_{D}^{-1} \widehat{C}_{e_n}(Y_i)
\]
and we estimate $\Gamma_I$ by the random operator
\[
\widehat{\Gamma}_{I,n}= \widehat{\Gamma}_{D}^{-1}\bigg( 2\widehat{\Gamma}_{e,n}+ \widehat{\Psi}_{e_n,D}-\Gamma_n\bigg).
\]
This random operator determines our   kernel FAVE approach for estimating the  EDR space. This estimation procedure is achieved  by considering   the space spanned by the  eigenvectors $ \widehat{\beta}_{1}, \widehat{\beta}_{2},\cdots, \widehat{\beta}_{K}$  of  $ \widehat{\Gamma}_{I,n},$ associated respectively with the $K$ largest eigenvalues  $ \widehat{\lambda}_{1},\cdots, \widehat{\lambda}_{K}$.
		\section{Asymptotics study of $ \widehat{\Gamma}_{I,n}$}
		In this section, we deal with asymptotics for $ \widehat{\Gamma}_{I,n}$. More precisely, we first establish its consistency as an estimator of $\Gamma_I$. Then we show the $\widehat{\beta}_k$'s are also consistent estimators of the $\beta_k$'s. We need the following assumptions:

\bigskip

			\noindent $(\mathscr{A}_5)$: $\Gamma$ is positive definite.

\bigskip

			\noindent$(\mathscr{A}_6)$: $f$,  $r$ and $R$ belong to $C^k$ ;

\bigskip
			\noindent$(\mathscr{A}_7)$: the kernel $K$ is of  order $k>2$ , has  compact support $[a,b]$, is  symmetric about zero and staisfies  $K\leq 1$, $\int_{a}^{b} |u|^kK(u)du<+\infty$, ;

\bigskip

			\noindent$(\mathscr{A}_8)$: there exist real numbers $d_1$, $d_2$ and $d_3$  such that  $\sup_{y\in\mathbb{R}}|f^{(k)}(y)|\leq d_1$, $\sup_{y\in\mathbb{R}}\Vert m^{(k)}(y)\Vert_{\mathcal{H}}\leq d_2$ and
			\noindent$\sup_{y\in\mathbb{R}}\Vert m^{(k)}(y)\Vert_{hs}\leq d_3$, where $\Vert\cdot\Vert_{hs}$ denotes the Hilbert-Scmidt norm of operators;

\bigskip

			\noindent $(\mathscr{A}_9)$:  $ h \sim n^{-c_1} $ and $e_n\sim n^{-c_2}$, where $c_1$ and $c_2$ are real numbers satisfying  $c_1>0$, $0<c_2<\frac{k-2}{4(k+1)}$ and  $\dfrac{c_2}{k}+\dfrac{1}{2k}<c_1<\dfrac{1}{4}-c_2$;

\bigskip

			\noindent$(\mathscr{A}_{10})$:  $\sqrt{n}\,\mathbb{E}\left[ \Vert R(Y)\Vert_{hs}^{2}\textbf{1}_{\{f(Y)<e_n\}}\right]$, $\sqrt{n}\,\mathbb{E}\left[ \Vert R(Y)\Vert_{hs}\, \Vert r(Y)\Vert_{H}^{2}\textbf{1}_{\{f(Y)<e_n\}}\right]$ and $\sqrt{n}\,\mathbb{E}\left[ \Vert R(Y)\Vert_{H}^{4}\textbf{1}_{\{f(Y)<e_n\}}\right]$ tends to $0$ as   $n\rightarrow +\infty$;

\bigskip

			\noindent$(\mathscr{A}_{11})$: the function $y\mapsto\mathbb{E}\left[\Vert X\Vert^2_\mathcal{H}|Y=y\right] $ is continuous.	
				\begin{remark}
					Zhu and Fang \cite{zhu} introduced $ \widehat{f}_{e_n}(y)=\max(\widehat{f}(y),e_n)$ to  overcome  technical diffiulties due to small values in the denominator of $ \widehat{r}(y)$. But this approach does not guarantee that we get
					a good estimator of $f$. Indeed,  if we take for example    $e_n=n^{-1/11}$, then until $n=2000$ we still have $e_n>1/2$ and, therefore,   $ \widehat{f}_{e_n}(y)=1/2$ for all $y\in\mathbb{R}$.   To overcome this later problem, Nkou and Nkiet\cite{nkou} propose to take $e_n=\min(a;n^{-c_2})$, where $a$ is a fixed strictly positive number. When $a$ is sufficiently
					small $ \widehat{f}_{e_n}(y)$ is a good estimator of $f$, because $\sup_{x\in \mathbb{R}}| \widehat{f}_{e_n}(y)- \widehat{f}(y)|\leq a$  and we still have $e_n\sim n^{-c_2}$ .  
				\end{remark}

			\noindent For $D\in\mathbb{N}^\ast$, we consider
\[
\Psi_D= \mathbb{E}\left[ Var(X|Y)\Gamma^{-1}_{D}Var(X|Y)\right]
\]
and denoting by   $t_D$  the minimum positive eigenvalue of $\Gamma_D$, we have:

			 \begin{Theo}\label{theo1}
				Under assumptions $(\mathscr{A}_{1})$ to $(\mathscr{A}_{3})$   and $(\mathscr{A}_{7})$ to $(\mathscr{A}_{11})$,  if we suppose that when $D\rightarrow +\infty$,  we have $\Vert \widehat{\Gamma}_{D}-\Gamma_{D}\Vert_\infty=o_p(t_{D}),$ then 
				\begin{align*}
				\Vert\Psi_D-\widehat{\Psi}_{e_n,D}\Vert_{hs}&=O_p\left( \dfrac{1}{\sqrt{n}}\right) +O_p\left( \dfrac{1}{t_{D}\sqrt{n}}\right)+O_p\left( \dfrac{1}{e_nt_{D}}\left( h^k+\dfrac{\sqrt{ \log(n)}}{h\sqrt{n}}\right) \right) \\
				&=O_p\left( \dfrac{1}{\sqrt{n}}\right)+O_p\left( \dfrac{1}{n^{\gamma}t_{D}}\right),
				\end{align*}
 where $\gamma$ is a constant satisfying $0<\gamma<1/4$.	
			 \end{Theo}
				\begin{remark}
				 This theorem gives an idea on the convergence rate of each component of $ \widehat{\Gamma}_{I,n}$ as we know the one of $ \widehat{\Gamma}_{e,n}$ from \cite{ferre}. We cannot reach the $\sqrt{n}$-convergence, because the rate of convergence will be penalized by the one of  $t_D$.The assumption    $\Vert \widehat{\Gamma}_{D}-\Gamma_{D}\Vert_\infty=o_p(t_{D})$ was also used in \cite{lian} for obtaining a similar result for the case of Functional SAVE. A justification of this asumption can be found in this paper.

				\end{remark}
\bigskip

\noindent In the following theorem consistency of   $ \widehat{\Gamma}_{I,n}$ is established under some conditions.
				 \begin{Theo}\label{theo2}
					Under the assumptions $(\mathscr{A}_{1})$ to $(\mathscr{A}_{11})$,  if we suppose that for some $0<\gamma<1/4,$ when $D\rightarrow +\infty$, $\Vert \widehat{\Gamma}_{D}-\Gamma_{D}\Vert_\infty=o_p(t_{D})$,  $1/(t_{D}\sqrt{n})\rightarrow 0,$ $1/(n^{\gamma}t_{D}^{2})\rightarrow 0,$ then  
					$\widehat{\Gamma}_{I,n}-\Gamma_I=o_p(1)$.
				 \end{Theo}
			\begin{remark}
			\noindent This result only gives the convergence  in probability of  $ \widehat{\Gamma}_{I,n}$  to $\Gamma_I$ without specifying	the rate.  For the functional SAVE, Lian and Li\cite{lian}  don't show the convergence of their estimator of $\Gamma_{I}$.
			\end{remark}
 \noindent Now, we deal with the $\widehat{\beta}_k$'s. For doing that, we  assume that $\beta_1,\beta_2,\cdots,\beta_K$ are the $K$ eigenvectors of $\Gamma_I$  associated   with the $K$ eigenvalues $\lambda_1,\cdots, \lambda_K$ respectiveley, and that  $\lambda_1>\lambda_2>\cdots>\lambda_K>0$.
			 \begin{Theo}\label{theo3}
						Under the assumptions $(\mathscr{A}_{1})$ to $(\mathscr{A}_{11})$,  if we suppose that for some $0<\gamma<1/4,$ when $D\rightarrow +\infty$, $\Vert \widehat{\Gamma}_{D}-\Gamma_{D}\Vert_\infty=o_p(t_{D})$, $1/(t_{D}\sqrt{n})\rightarrow 0,$ $1/(n^{\gamma}t_{D}^{5/2})\rightarrow 0,$ then $\Vert \widehat{\beta}_{j}-\beta_j\Vert_{\mathcal{H}}=o_p(1)$ for $j=1,2,\cdots ,K$.
			 \end{Theo}
			\begin{remark}
			This result is similar to that of FSIR obtained by Ferr\'e  and Yao\cite{ferre}. It is an extension to the functional case   of a property of  the kernel method for sliced average variance estimation developped by Zhu and Zhu\cite{zhu2} in a multivariate context. 
			\end{remark}
\section{Simulation study}
				In this section, we use  simulations  to illustrate the kernel  FAVE method and to compare it with existing methods. In all the examples, the predictor $X$ is a standard brownian motion on $[0,1]$, observed on a grid of  $p=100$ equally spaced points.  Two models are considered:

\bigskip

\noindent\textbf{Model 1:} $ Y=\sin(\pi<\beta_1,X>_{\mathcal{H}}/2)+<\beta_2,X>_{\mathcal{H}}^{5}+\varepsilon$,  where $\beta_1(t)=(2t-1)^3+1$, $\beta_2(t)=\cos(\pi(2t-1))+1$ and $\varepsilon\sim N(0,0.1^2)$.

\bigskip

\noindent\textbf{Model 2:} $ Y=50<\beta_1,X>_{\mathcal{H}}^{2}+<\beta_2,X>_{\mathcal{H}}^{2}+\varepsilon$,  where $\beta_1(t)=4t^2$, $\beta_2(t)=\sin(5\pi t/2)$ and $\varepsilon\sim N(0,0.1^2)$.

\bigskip

\noindent We set $n=100$ and we consider both functional PCA and quadratic  B-spline basis for computing $\widehat{\Pi}_D$. For the B-spline basis, the knots are chosen to be equally spaced on $[0,1]$.  Various dimensions $D$ are used, $D=4,5,6,7,8$.  The bandwidth $h$ is selected by the cross-validation. 

\bigskip
				 
				  \noindent  The  plots of $\hat{\beta}_1$ and $\hat{\beta}_2$, obtained from kernel FAVE, together with that of  $\beta_1$ and $\beta_2$ are given in Figure 1-2 for  Model 1, and in Figure 3-4 for Model 2. They reveal very good estimations. In order to verify if the prior projection space is well estimates by our FAVE method, we plot in Figure 5 to 8, the index $<\beta_j,X>_\mathcal{H}$  versus $<\hat{\beta}_j,X>_\mathcal{H}$   for $j=1,2$ and for Model 1-2.  These scatter plots reveal a strong correlation in both cases. All the previous plot are made using $D=4$.\\
				  
				  \noindent In order to compare our method to the FSIR and FSAVE methods, we use various dimensions $D=4$ to $D=8$.  FSAVE is performed with number of slices   $H=10$. The distance between the true  EDR space and its estimation is computed via   $\mathscr{E}=\Vert P-\hat{P}\Vert_{hs}$, where   $P$  (resp. $\widehat{P}$) denotes the projector onto the space spanned by $\beta_1$ and $\beta_2$ (resp.  $\widehat{\beta}_1$ and $\widehat{\beta}_2$). We use $m=100$ simulated datasets in each scenarios to get the boxplot of $\mathscr{E}$. In the left hand of each, figure from 9 to 14, the boxplots are built   by using functional PCA basis expansion, whereas the ones in the rigth hand are  based on   B-spline basis functions.  The boxplot results  related to the B-spline projections are almost better than the ones from functional PCA. For Model 1 the three methods perform similarly as showed by  the boxplots, but   in the case of Model 2 FSIR does not work  as well as FSAVE and FAVE. As a general observation the three methods are sensitives to the choice of $D$. Therefore, methods are needed for  chosing  of $D$ and  will perfect the pratical  use  of   the FSIR, FAVE and FSAVE methods.
  
\begin{figure*}
	\begin{minipage}[h]{.45\linewidth}
		\centering
		\includegraphics[width=1.0\linewidth]{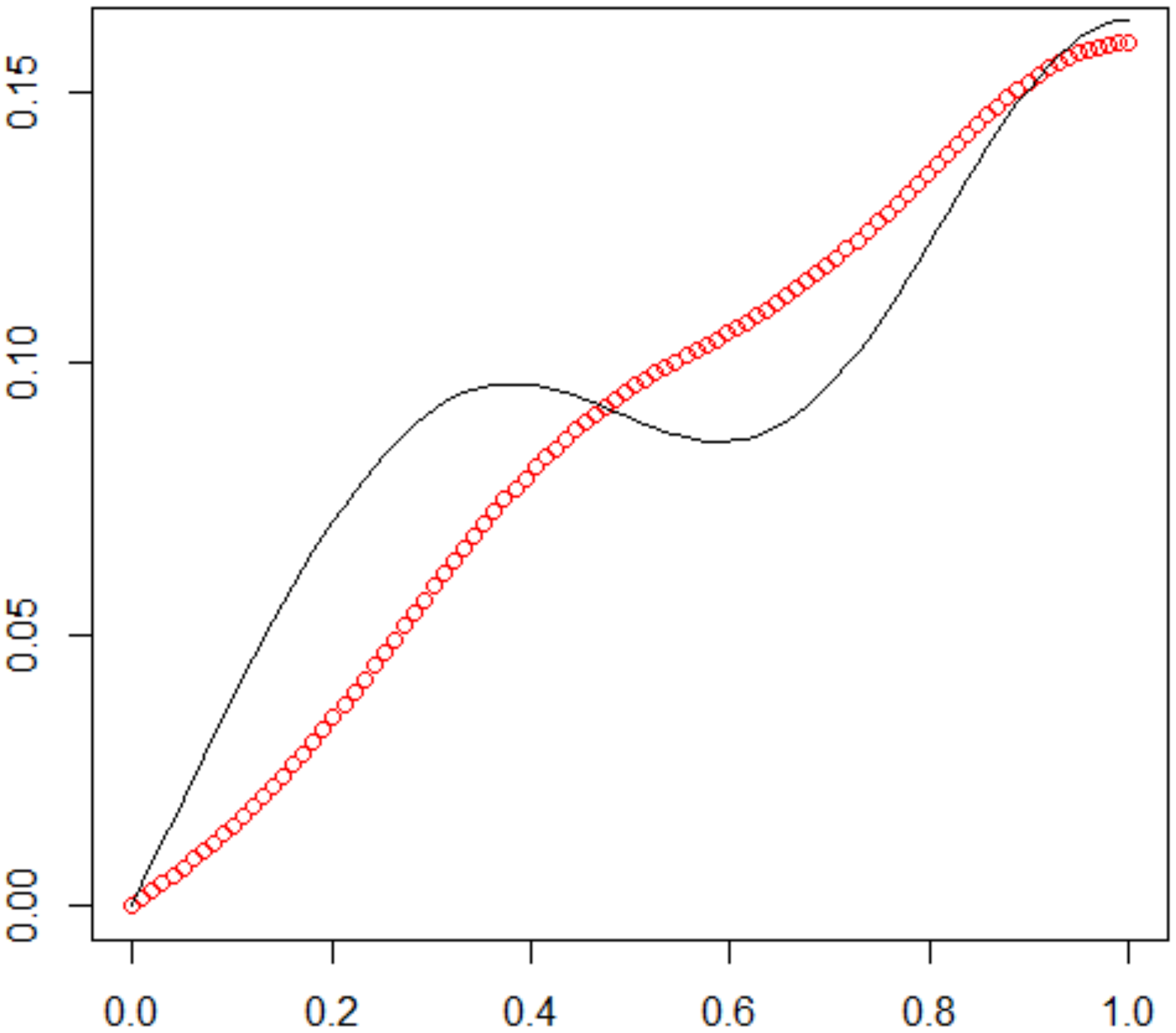}
		\caption{Plot of $\widehat{\beta}_1$ (red) and $\beta_1$ (black) for Model 1.}
	\end{minipage}
	\label{fig:M1_first_EDR_direction_FAVE}
	\hfill
		\begin{minipage}[h]{.45\linewidth}
			\centering
			\includegraphics[width=1.0\linewidth]{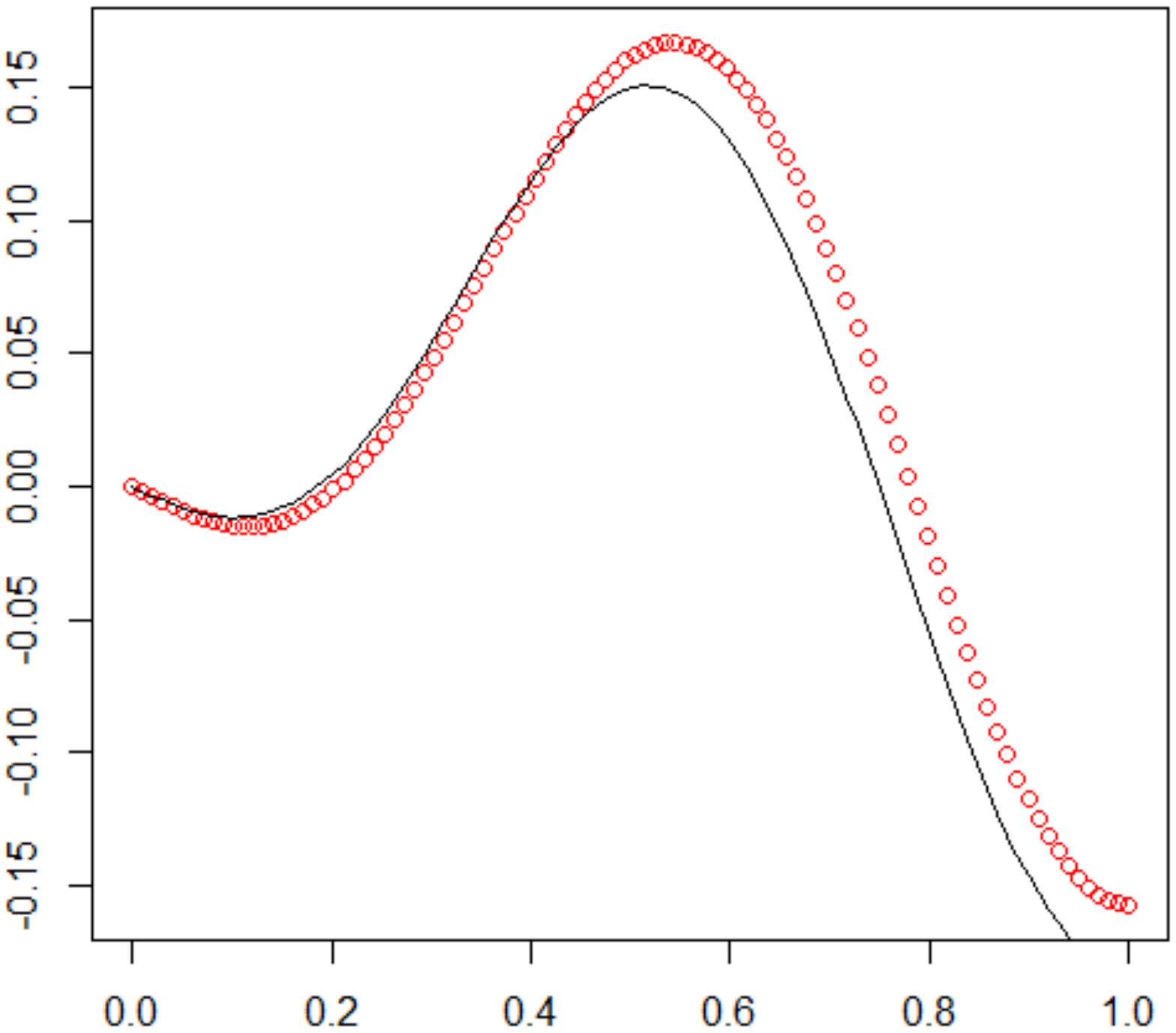}
			\caption{Plot of $\widehat{\beta}_2$ (red) and $\beta_2$ (black) for Model 1.}
		\end{minipage}
		\label{fig:M1_second_EDR_direction_FAVE}
		\hfill
	\begin{minipage}[h]{.45\linewidth}
\centering
\includegraphics[width=1.0\linewidth]{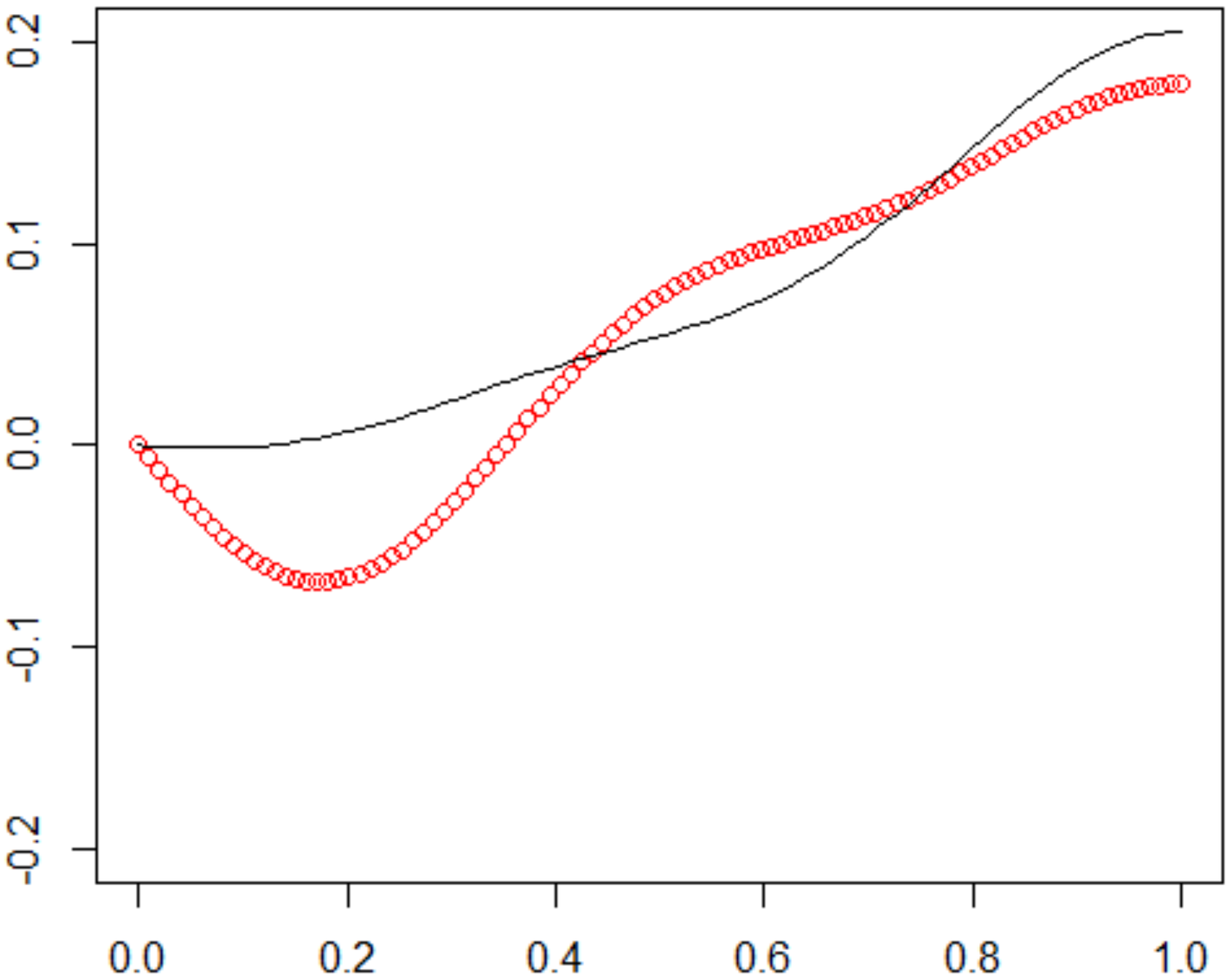}
\caption{Plot of $\widehat{\beta}_1$ (red) and $\beta_1$ (black) for Model 2.}
	\end{minipage}
\label{fig:M2_first_EDR_direction_FAVE}
\hfill
	\begin{minipage}[h]{.45\linewidth}
		\centering
		\includegraphics[width=1.0\linewidth]{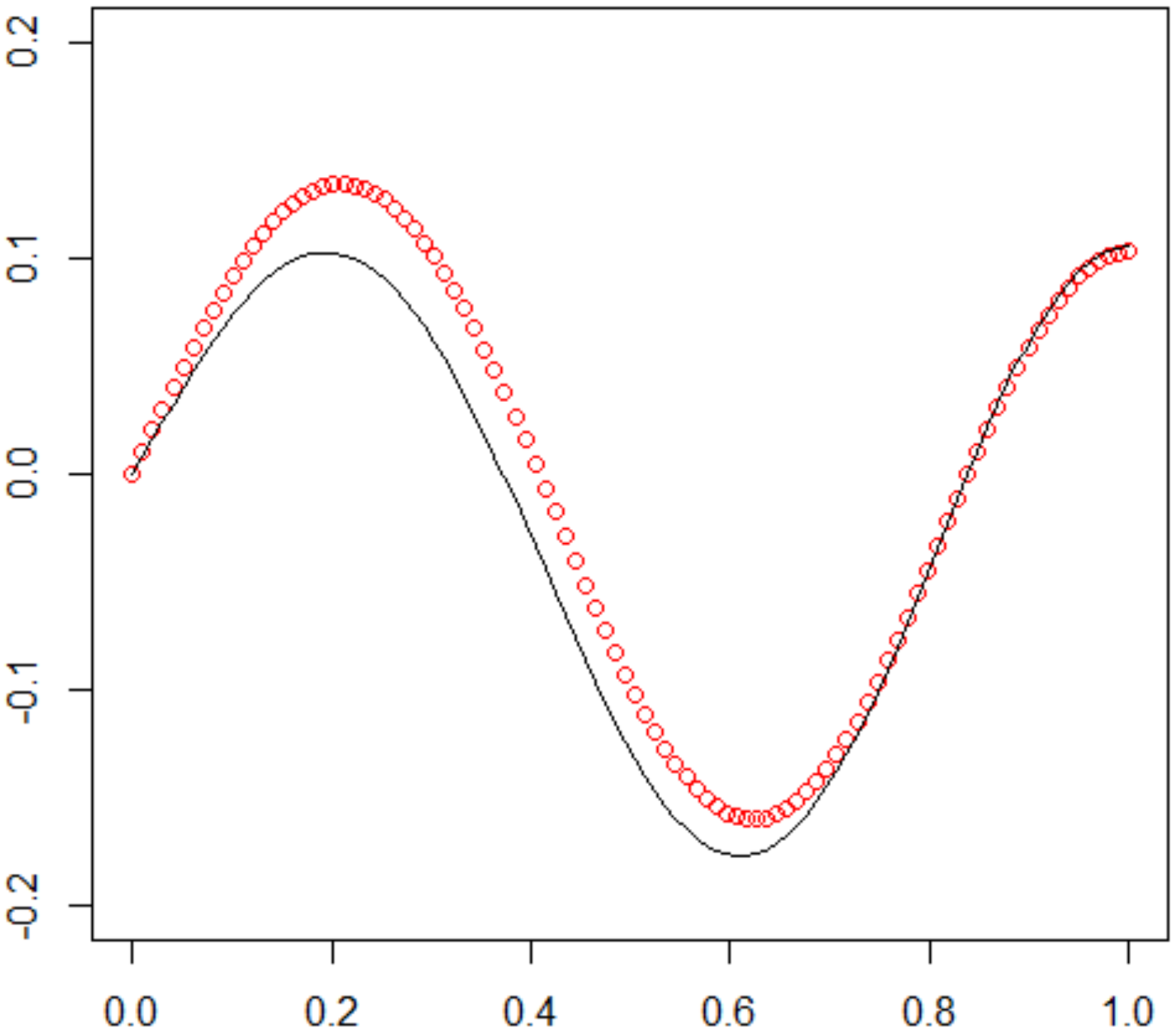}
		\caption{Plot of $\widehat{\beta}_2$ (red) and $\beta_2$ (black) for Model 2.}
	\end{minipage}
	\label{fig:M2_second_EDR_direction_FAVE}
\end{figure*}

\begin{figure}
	\begin{minipage}[h]{.45\linewidth}
		\centering
		\includegraphics[width=1.0\linewidth]{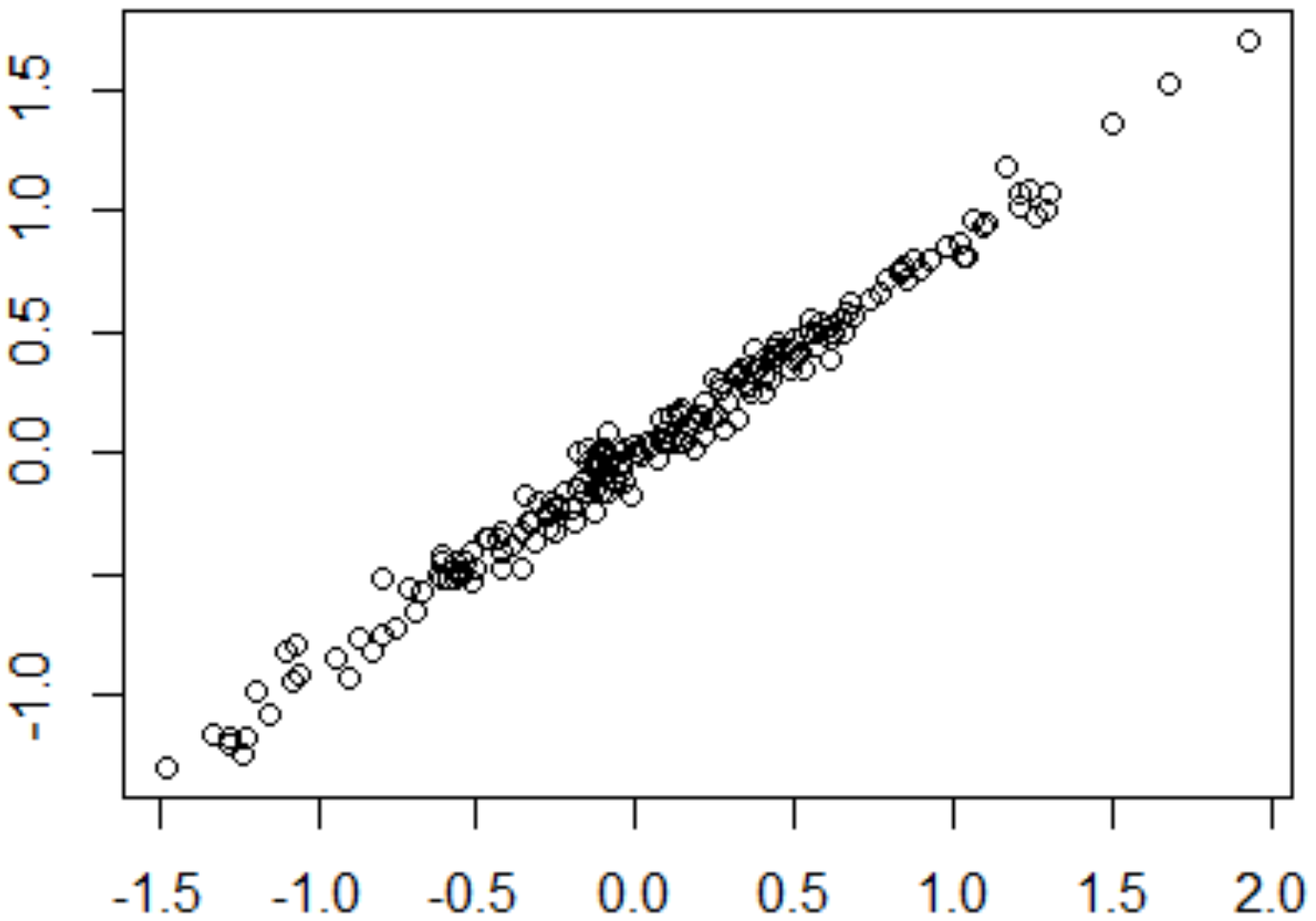}
		\caption{Plot of $<\widehat{\beta}_1,X>_\mathcal{H}$  versus $<\beta_1,X>_\mathcal{H}$  for Model 1.}
	\end{minipage}
	\label{fig:M1_first_indexe}
		\hfill
	\begin{minipage}[h]{.45\linewidth}
\centering
\includegraphics[width=1.0\linewidth]{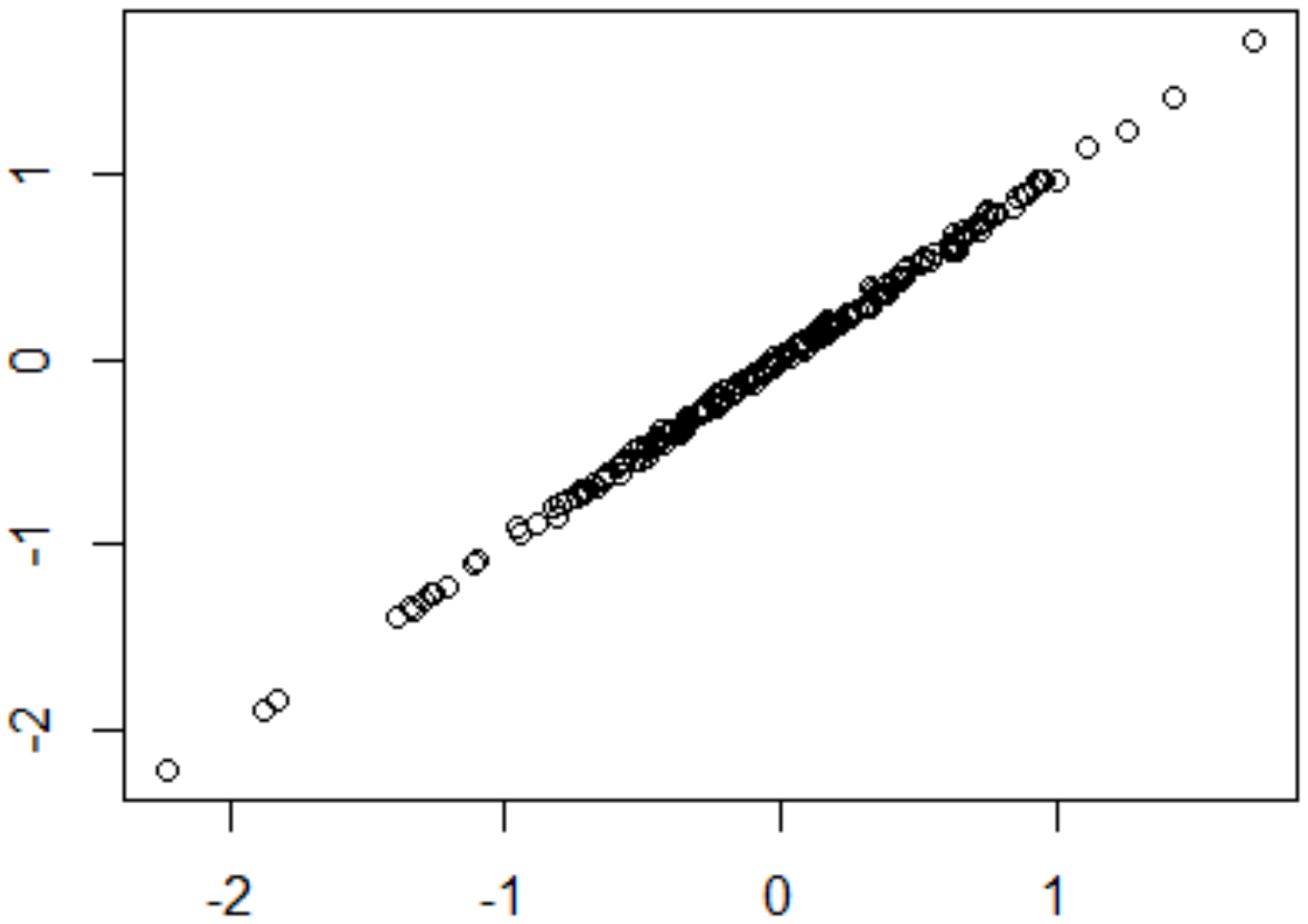}
\caption{Plot of $<\widehat{\beta}_2,X>_\mathcal{H}$  versus  $<\beta_2,X>_\mathcal{H}$  for Model 1.}
\end{minipage}
\label{fig:M1_second_indexe}
	\hfill
	\begin{minipage}[h]{.45\linewidth}
		\centering
		\includegraphics[width=1.0\linewidth]{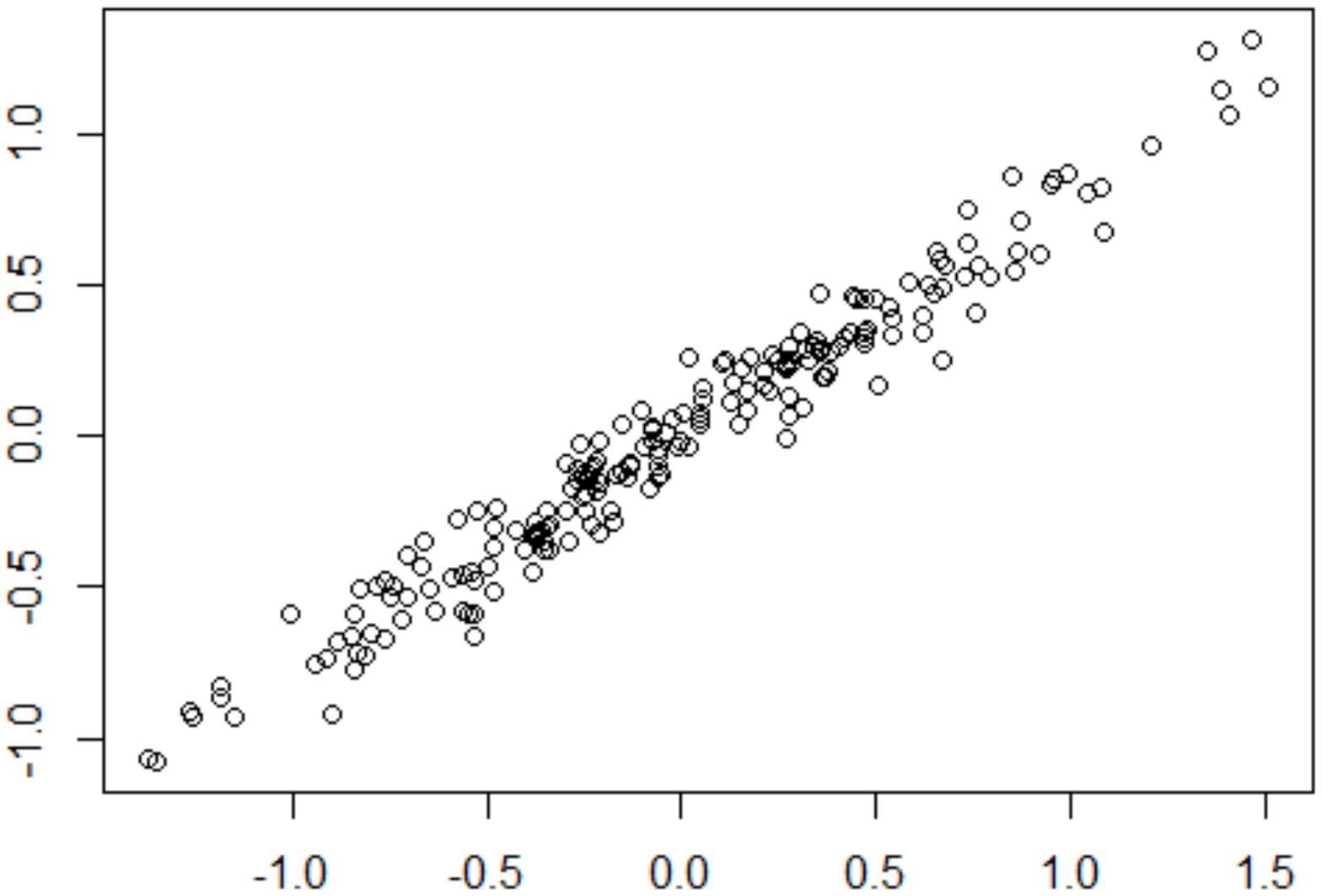}
		\caption{Plot of $<\widehat{\beta}_1,X>_\mathcal{H}$  versus $<\beta_1,X>_\mathcal{H}$  for Model 2.}	
	\end{minipage}
	\label{fig:M2_first_indexes}
		\hfill
	\begin{minipage}[h]{.45\linewidth}
		\centering
		\includegraphics[width=1.0\linewidth]{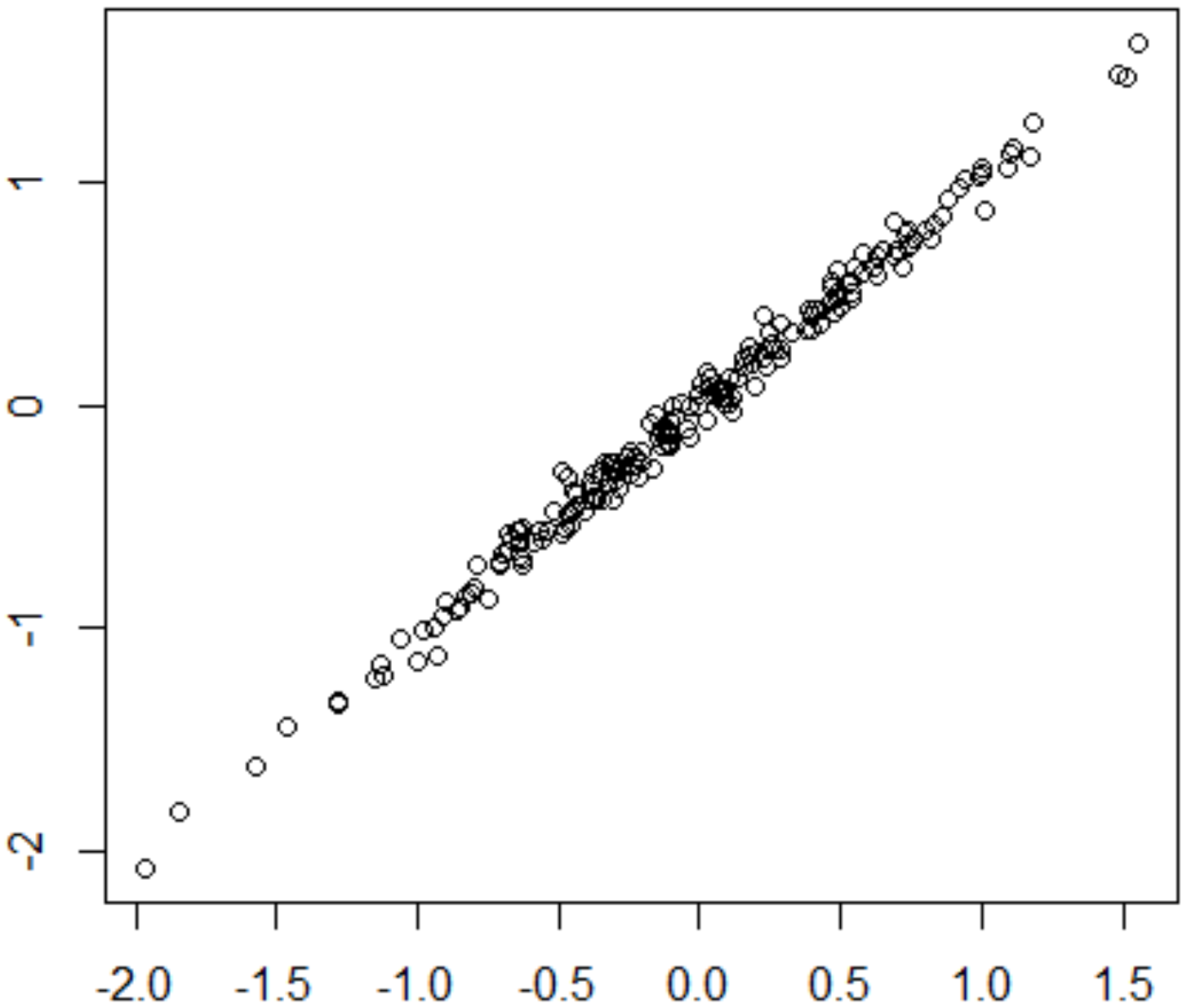}
		\caption{Plot of $<\widehat{\beta}_2,X>_\mathcal{H}$  versus  $<\beta_2,X>_\mathcal{H}$  for Model 2.}	
	\end{minipage}
	\label{fig:M2_second_indexes}
\end{figure}

\begin{figure}
	
	\begin{minipage}[h]{.55\linewidth}
		\centering
		\includegraphics[width=1.0\linewidth]{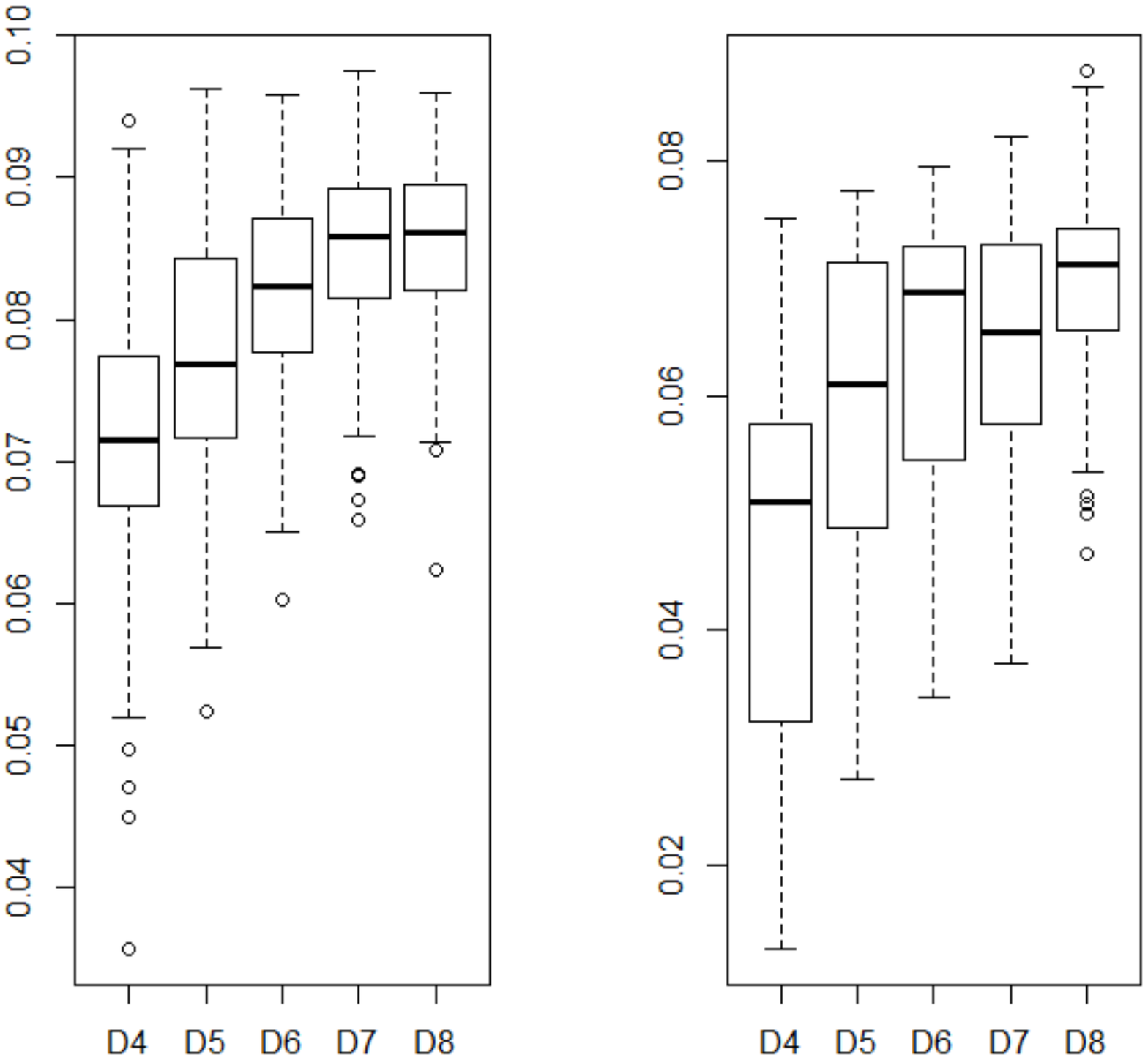}
		\caption{Boxplots showing $\Vert P-\hat{P}\Vert_{hs}$ for Model 1, from FSAVE with $H=10$.}
	\end{minipage}
	\label{fig:M1_FSAVE_H=10}
	%\hfill
	
\end{figure}
\begin{figure}
		\begin{minipage}[h]{.55\linewidth}
\centering
\includegraphics[width=1.0\linewidth]{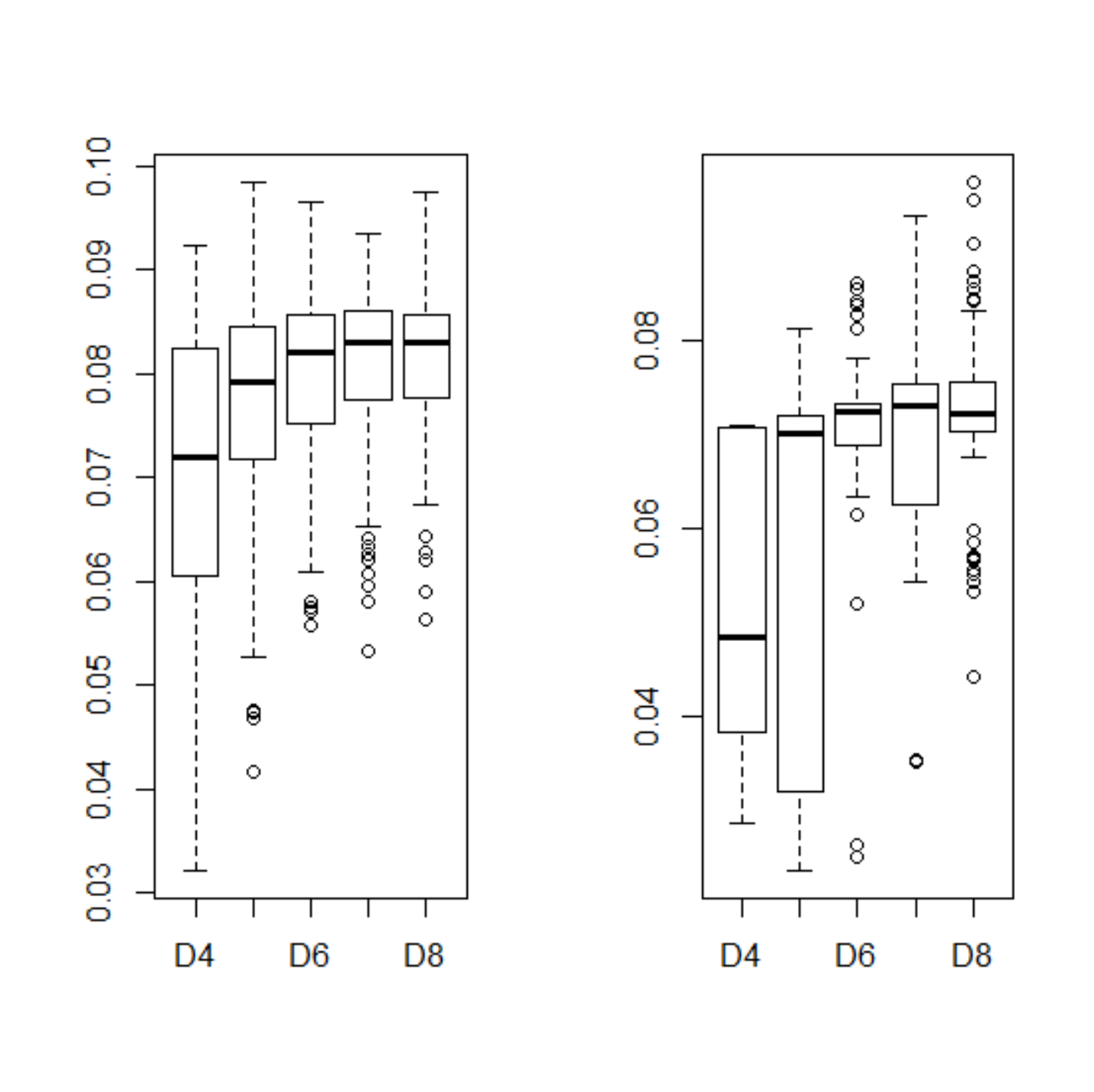}
\caption{Boxplots showing $\Vert P-\hat{P}\Vert_{hs}$ for Model 1, using  kernel FAVE.}
	\end{minipage}
\label{fig:M1_boxplot_FAVE}
\hfill
	\begin{minipage}[h]{.55\linewidth}
		\centering
		\includegraphics[width=1.0\linewidth]{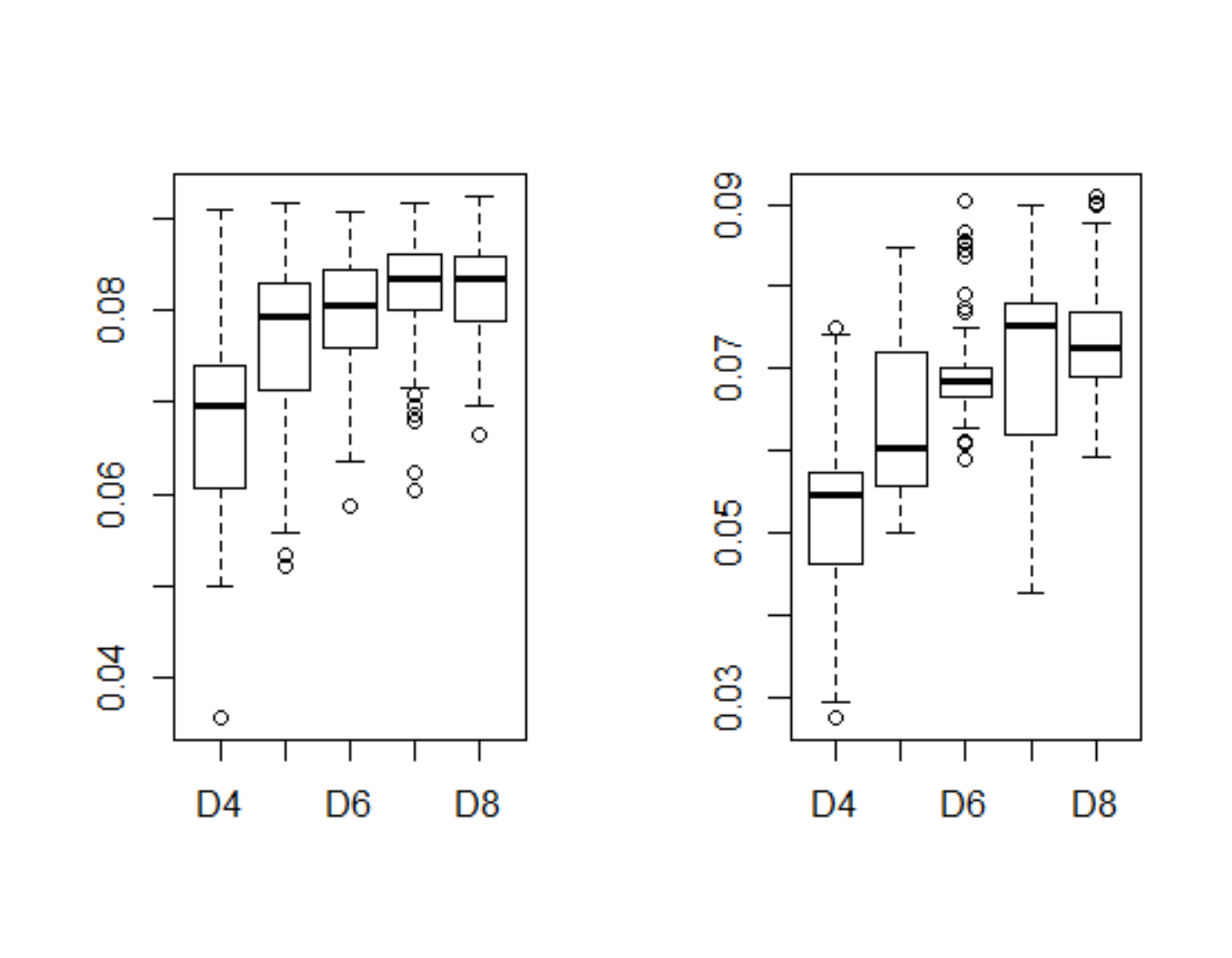}
		\caption{Boxplots showing $\Vert P-\hat{P}\Vert_{hs}$ for Model 1, using  kernel FSIR.}
	\end{minipage}
	\label{fig:M1_boxplot_FIR}
\end{figure}

\begin{figure}
	\begin{minipage}[h]{.55\linewidth}
		\centering
		\includegraphics[width=1.0\linewidth]{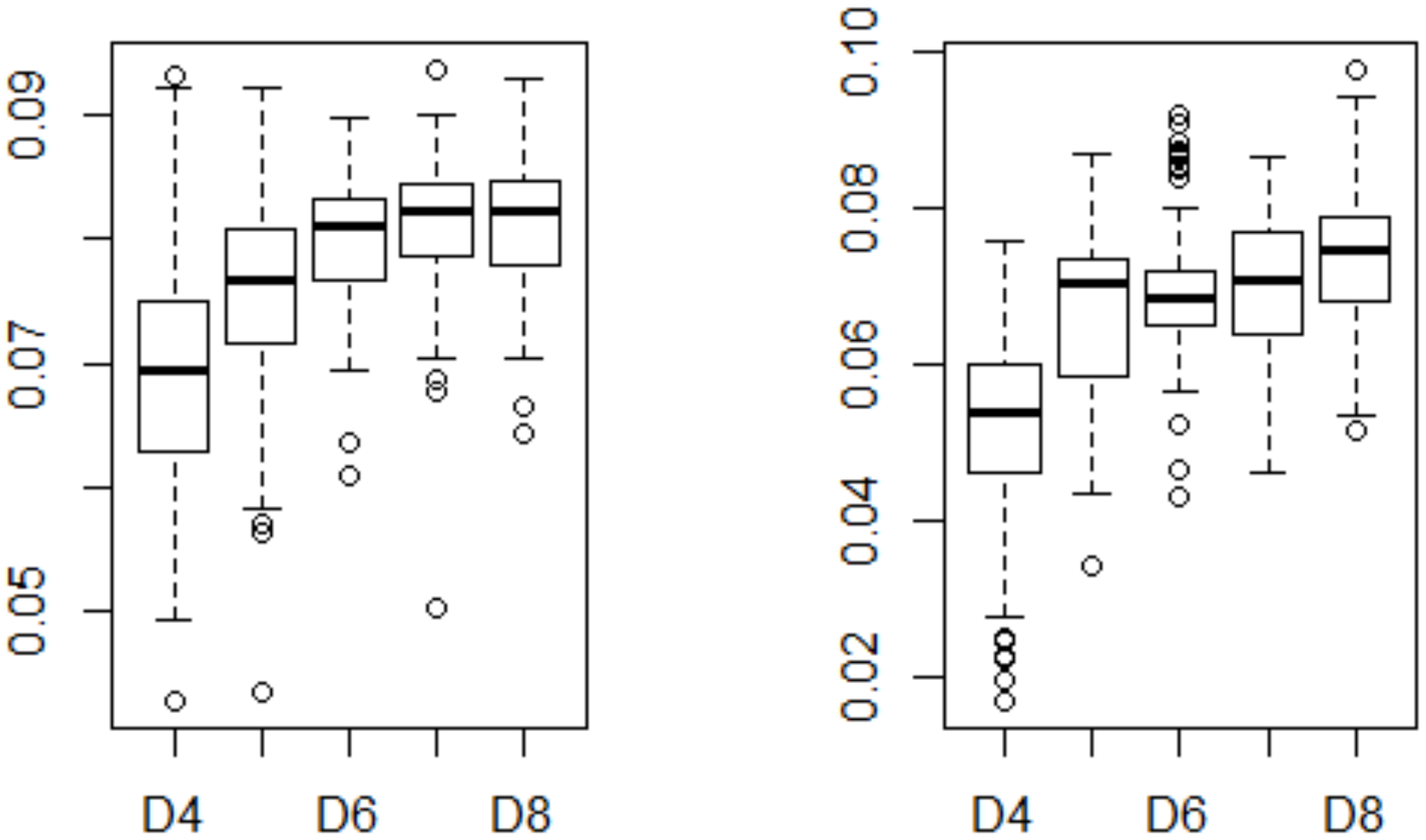}
		\caption{Boxplots showing $\Vert P-\hat{P}\Vert_{hs}$ for Model 2, from FSAVE with $H=10$}
	\end{minipage}
	\label{fig:M2_FSAVE_H=10}
		%\hfill
\end{figure}
\begin{figure}
		\begin{minipage}[h]{.55\linewidth}
\centering
\includegraphics[width=1.0\linewidth]{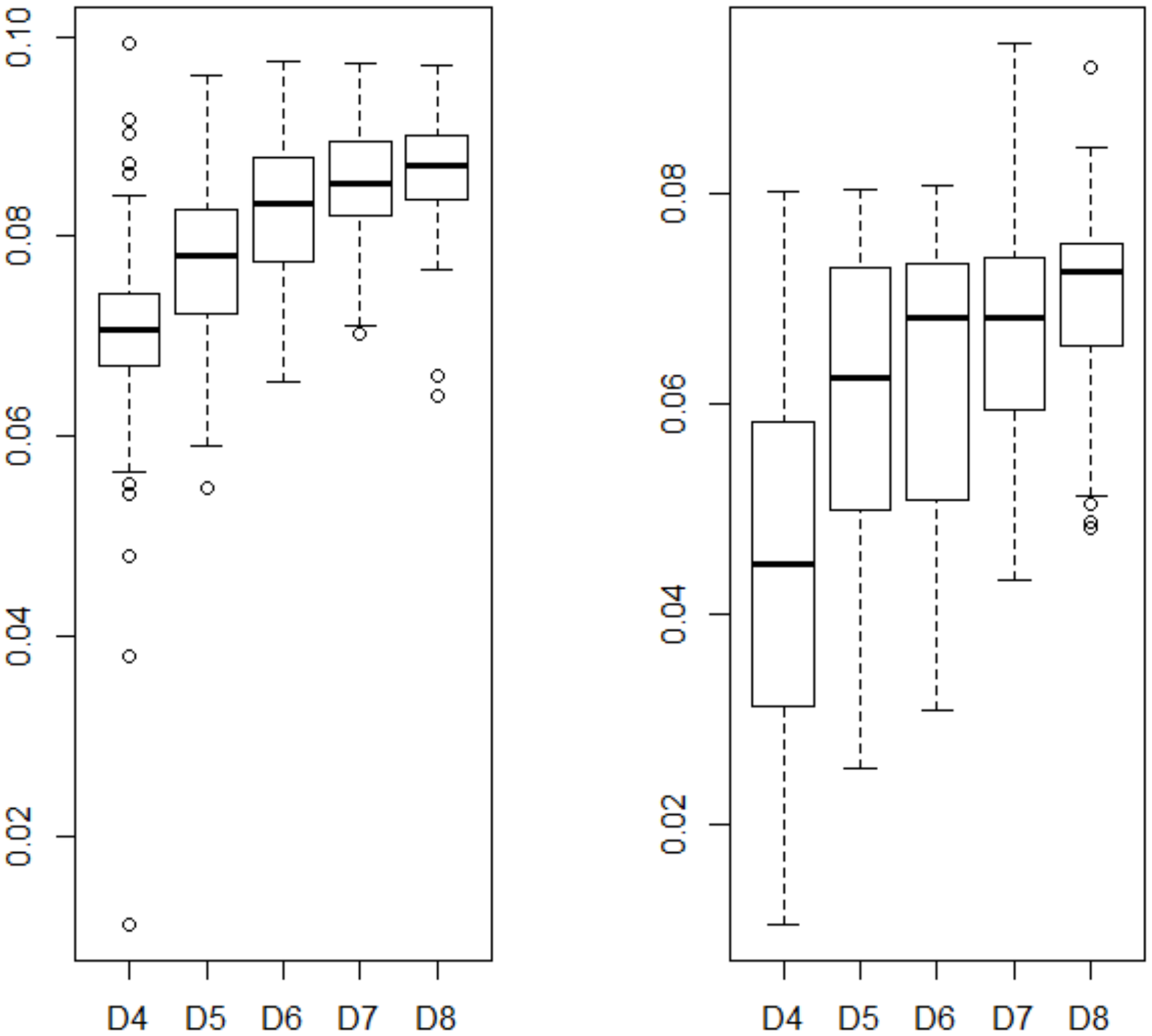}
\caption{Boxplots showing $\Vert P-\hat{P}\Vert_{hs}$ for Model 2, using  kernel FAVE.}
	\end{minipage}
\label{fig:M2_boxplot_FAVE}
	\begin{minipage}[h]{.55\linewidth}
		\centering
		\includegraphics[width=1.0\linewidth]{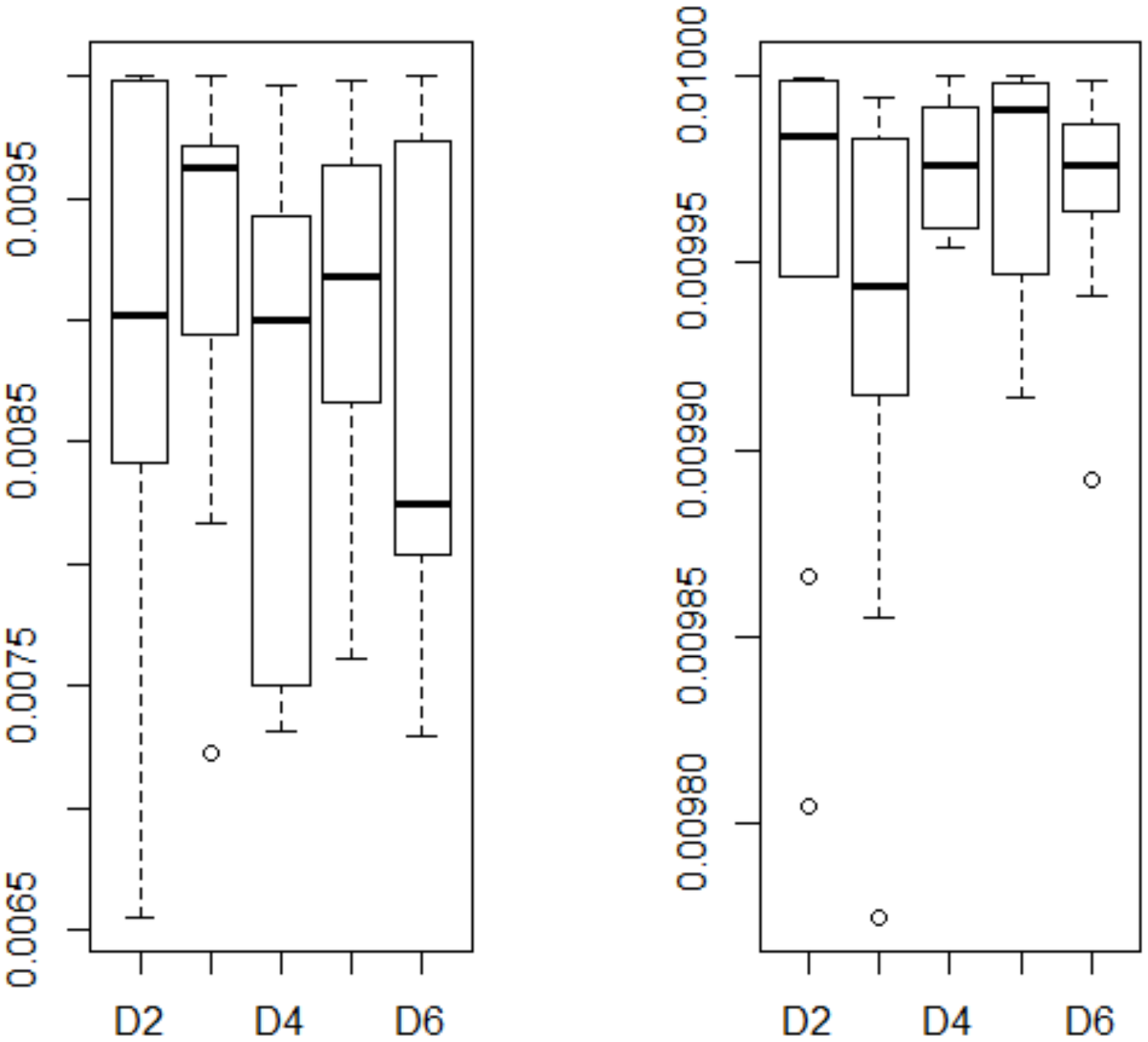}
		\caption{Boxplots showing $\Vert P-\hat{P}\Vert_{hs}$ for Model 1, using  kernel FSIR.}
	\end{minipage}
	\label{fig:M2_Boxplot_FIR}
\end{figure}

\section{Proofs}
\subsection{Preliminary results}
In this section we will give some lemmas necessary to get the proofs of the previous Theorems. 

\begin{lemma}\label{lem2}
	Under assumption $(\mathscr{A}_6)$ to  $(\mathscr{A}_9)$, we have
	\[
	\sup_{y\in \mathbb{R}}\Vert \widehat{M}(y)-M(y)\Vert_{hs}=O_p\left( h^k+\dfrac{\sqrt{ \log(n)}}{h\sqrt{n}}\right). 
	\]
\end{lemma}
\noindent\textit{Proof.}
It is easy to check that for all  $y\in \mathbb{R}$, one has     
\[
\mathbb{E}\left[ \widehat{M}(y) \right] =\dfrac{M*K_h(y)}{h}.
\] 
Then
\begin{align*}
\mathbb{E}\left[ \widehat{M}(y) \right]-M(y) &= \dfrac{1}{h}\int_{\mathbb{R}}R(v)K_h(v-y)f(v)dv-M(y)\\
&=\dfrac{1}{h}\int_{\mathbb{R}}\left[ M(v)-M(y)\right] K_h(v-y)dv\\
	&=\int_{\mathbb{R}}\left[ M(y+hw)-M(y)\right] K(w)dw \\
          &= \int_{\mathbb{R}}\left[\sum_{j=1}^{k-1}\dfrac{(wh)^j}{j!}M^{(j)}(y)+\dfrac{(wh)^k}{k!}M^{(k)}(y+\theta hw) \right] K(w)dw\\
          &=\dfrac{h^k}{k!}\int_{a}^{b}w^kM^{(k)}(y+\theta hw)  K(w)dw.
\end{align*}
Hence
\begin{align*}
	\Vert\mathbb{E}\left[ \widehat{M}(y) \right]-M(y)\Vert_{hs}&\leq\dfrac{h^k}{k!}\sup_{y\in I}\Vert M^{(k)}(y)\Vert_{hs} \int_{a}^{b}|w|^k K(w)dw=Ch^k ,
\end{align*}
that is $\sup_{y\in I}\Vert\dfrac{M*K_h(y)}{h}-M(y)\Vert_{hs}=O(h^k)$. We deduce that
	\begin{align*}
	\sup_{y\in \mathbb{R}}\Vert \widehat{M}(y)-M(y)\Vert_{hs}&\leq \sup_{y\in \mathbb{R}}\Vert	 \widehat{M}(y)-\mathbb{E}\left[ \widehat{M}(y) \right]\Vert_{hs}+\sup_{y\in \mathbb{R}}\Vert\mathbb{E}\left[ \widehat{M}(y) \right]-M(y)\Vert_{hs}\\
	&=D_1+D2.
	\end{align*}
	Let $\varepsilon>0$ and $(a_n)_{n\in \mathbb{N}},$ a sequence of non-negative reals numbers converging to $+\infty$. We have :
	\begin{align*}
	P\left( D_1>\varepsilon\right) &=P\left( \sup_{y\in \mathbb{R}}\Vert	 \widehat{M}(y)-\mathbb{E}\left[ \widehat{M}(y) \right]\Vert_{hs}>\varepsilon\right)\\
	&\leq P\left( \sup_{y\in \mathbb{R}}\Vert	 \widehat{M}(y)-\mathbb{E}\left[ \widehat{M}(y) \right]\Vert_{hs}>\varepsilon;\Vert X\otimes X\Vert_{hs}\leq a_n\right)\\
	& +P\left(\sup_{y\in \mathbb{R}}\Vert	 \widehat{M}(y)-\mathbb{E}\left[ \widehat{M}(y) \right]\Vert_{hs}>\varepsilon;\Vert X\otimes X\Vert_{hs}> a_n\right)\\
	&\leq P\left( \sup_{y\in \mathbb{R}}\Vert	 \widehat{M}(y)-\mathbb{E}\left[ \widehat{M}(y) \right]\Vert_{hs}>\varepsilon;\Vert X\otimes X\Vert_{hs}\leq a_n\right)\\
	&+ P\left(\Vert X\otimes X\Vert_{hs}> a_n\right).
	\end{align*}
	Since $K\leq 1$, we have for all $ y\in \mathbb{R}$
	\begin{align*}
	\Vert	 \widehat{M}(y)-\mathbb{E}\left[ \widehat{M}(y) \right]\Vert_{hs}&=\Vert\dfrac{1}{nh}\sum\limits_{i=1}^{n}\left[ X_i\otimes X_iK\left(\dfrac{Y_i-y}{h} \right) -\mathbb{E}[ X_i\otimes X_iK\left(\dfrac{Y_i-y}{h}\right) ] \right] \Vert_{hs}\\
	&\leq \dfrac{1}{nh}\sum\limits_{i=1}^{n}\{\Vert X_i\otimes X_i\Vert_{hs} +\mathbb{E}[ \Vert X_i\otimes X_i \Vert_{hs}] \}.
\end{align*}
Thus
\begin{align*}
	\sup_{y\in \mathbb{R}}\Vert	 \widehat{M}(y)-\mathbb{E}\left[ \widehat{M}(y) \right]\Vert_{hs}&\leq\dfrac{1}{nh}\sum\limits_{i=1}^{n}\{\Vert X_i\otimes X_i\Vert_{hs} +\mathbb{E}[ \Vert X_i\otimes X_i \Vert_{hs}] \}
\end{align*}
and
	\begin{align*}
	& P\left( \sup_{y\in \mathbb{R}}\Vert	 \widehat{M}(y)-\mathbb{E}\left[ \widehat{M}(y) \right]\Vert_{hs}>\varepsilon;\Vert X\otimes X\Vert_{hs}\leq a_n\right)\\
	&\leq P\left( \dfrac{1}{nh}\sum\limits_{i=1}^{n}\left( \Vert X_i\otimes X_i\Vert_{hs} +\mathbb{E}[ \Vert X_i\otimes X_i\Vert_{hs}] \right)>\varepsilon;\Vert X\otimes X\Vert_{hs}\leq a_n\right)\\
	&\leq P\left( \dfrac{1}{nh}\sum\limits_{i=1}^{n}\left( \Vert X_i\otimes X_i\Vert_{hs} +\mathbb{E}[ \Vert X_i\otimes X_i\Vert_{hs}] \right)\textbf{1}_{\{\Vert X\otimes X\Vert_{hs}\leq a_n\}}>\varepsilon\right).
	\end{align*}
	However,  for any $ i\in\{1,\cdots,n\}$, one has
	\begin{align*}
	\dfrac{1}{h}\left( \Vert X_i\otimes X_i\Vert_{hs} +\mathbb{E}[ \Vert X_i\otimes X_i\Vert_{hs}] \right)\textbf{1}_{\{\Vert X\otimes X\Vert_{hs}\leq a_n\}}&\leq \dfrac{2a_n}{h}.
	\end{align*}
	Then using Bernstein inequality, we obtain
	\begin{align*}
	P\left( \sup_{y\in \mathbb{R}}\Vert	 \widehat{M}(y)-\mathbb{E}\left[ \widehat{M}(y) \right]\Vert_{hs}>\varepsilon;\Vert X\otimes X\Vert_{hs}\leq a_n\right)&\leq 2exp\left( -\dfrac{n\varepsilon^2 h^2 }{16a_{n}^{2}} \right) 
	\end{align*} 
from what we deduce
	\begin{align*}
	P\left(D_2>\varepsilon\right) &\leq P\left(\Vert X\otimes X\Vert_{hs}> a_n\right)+2\exp\left( -\dfrac{n\varepsilon^2 h^2 }{16a_{n}^{2}} \right)\\
	&\leq \dfrac{\mathbb{E}\left(\Vert X\otimes X\Vert_{hs}^{2}\right)}{a_{n}^{2}}+2\exp\left(- \dfrac{n\varepsilon^2 h^2 }{16a_{n}^{2}} \right)\\
	&= \dfrac{\mathbb{E}\left(\Vert X\Vert_{H}^{4}\right)}{a_{n}^{2}}+2\exp\left( -\dfrac{n\varepsilon^2 h^2 }{16a_{n}^{2}} \right).
\end{align*}
	Taking  $\varepsilon=\dfrac{\varepsilon_0}{h}\sqrt{\dfrac{\log(n)}{n}}$, where $\varepsilon >0$, and $ a_n=\left(\log (n)\right) ^{1/4}$, we have
\begin{align*}
	P\left(D_2>\dfrac{\varepsilon_0}{h}\sqrt{\dfrac{\log(n)}{n}}\right) &\leq  \dfrac{\mathbb{E}\left(\Vert X\Vert_{H}^{4}\right)}{\left( \log(n)\right) ^{1/2}}+2\exp\left( -\dfrac{\varepsilon_0(\log(n))^{1/2}  }{16}\right)
\end{align*}
and since 
\begin{equation*}
\lim\limits_{n\longrightarrow +\infty}\left( \dfrac{\mathbb{E}\left(\Vert X\Vert_{H}^{4}\right)}{\left(  \log(n)\right) ^{1/2}}+2exp\left( -\dfrac{\varepsilon_0( \log(n))^{1/2}  }{16} \right) \right) =0
\end{equation*}
we conclude that $	D_2=O_p\left( h^{-1}n^{-1/2}(\log(n))^{1/2} \right)$ and, consequently, that 
\[
	\sup_{y\in \mathbb{R}}\Vert \widehat{M}(y)-M(y)\Vert_{hs}=O_p\left( h^k+\dfrac{\sqrt{ \log(n)}}{h\sqrt{n}}\right).
\]
\hfill $\Box$

	\begin{lemma}\label{lem3}
We have:
\[
		\left| \dfrac{ \widehat{f}_{e_n}(Y_j)- \widehat{f}(Y_j)}{f_{e_n}(Y_j)}\right|\leq 2\left[ \textbf{1}_{\{f(Y_j)< 2e_n\}}+\dfrac{\left( \sup_{y\in \mathbb{R}}| \widehat{f}(y)-f(y)|\right) ^2}{e_{n}^{2}}\right].
\]
	\end{lemma}
\noindent\textit{Proof.}
Since
	\begin{align*}
	| \widehat{f}_{e_n}(Y_j)- \widehat{f}(Y_j)|&=|e_n\textbf{1}_{\{ \widehat{f}(Y_j)<e_n\}}+ \widehat{f}(Y_j)\textbf{1}_{\{ \widehat{f}(Y_j)\geq e_n\}}- \widehat{f}(Y_j)|\\
	&=|e_n\textbf{1}_{\{ \widehat{f}(Y_j)<e_n\}}- \widehat{f}(Y_j)\textbf{1}_{\{ \widehat{f}(Y_j)< e_n\}}|\\
	&\leq\left( e_n+ \widehat{f}(Y_j)\right)\textbf{1}_{\{ \widehat{f}(Y_j)<e_n\}}
\end{align*}
we obtain
\begin{align*}
	\left| \dfrac{ \widehat{f}_{e_n}(Y_j)- \widehat{f}(Y_j)}{f_{e_n}(Y_j)}\right| &\leq \left( \dfrac{e_n}{f_{e_n}(Y_j)}+\dfrac{ \widehat{f}(Y_j)}{f_{e_n}(Y_j)}\right) \textbf{1}_{\{ \widehat{f}(Y_j)< e_n\}}\leq 2\,  \textbf{1}_{\{ \widehat{f}(Y_j)< e_n\}}.
	\end{align*}
	It is easy to check that 
\[
\textbf{1}_{\{ \widehat{f}_{e_n}(Y_j)< e_n\}}\leq \textbf{1}_{\{f(Y_j)< 2e_n\}}+\dfrac{\left( \sup_{y\in \mathbb{R}}| \widehat{f}(y)-f(y)|\right) ^2}{e_{n}^{2}}.
\]
Thus
\[
	\left| \dfrac{ \widehat{f}_{e_n}(Y_j)- \widehat{f}(Y_j)}{f_{e_n}(Y_j)}\right|\leq 2\left[ \textbf{1}_{\{f(Y_j)< 2e_n\}}+\dfrac{\left( \sup_{y\in \mathbb{R}}| \widehat{f}(y)-f(y)|\right) ^2}{e_{n}^{2}}\right].
\]
\hfill $\Box$

	\begin{lemma}\label{lem4}
		Under the assumption ($\mathscr{A}_3)$,  if we suppose $\Vert \widehat{\Gamma}_{D}-\Gamma_{D}\Vert_\infty=o_p(t_{D})$, we have
		$$ A_{2n}=\dfrac{1}{n}\sum_{j=1}^{n}C(Y_j)\Gamma^{-1}_{D}C(Y_j)-\dfrac{1}{n}\sum_{j=1}^{n}C(Y_j) \widehat{\Gamma}^{-1}_{D}C(Y_j)=O_p\left( \dfrac{1}{t_{D}\sqrt{n}}\right).$$
	\end{lemma}
\noindent\textit{Proof.}

	\begin{align*}
	\Vert A_{2n}\Vert_{hs}&= \Vert\dfrac{1}{n}\sum_{j=1}^{n}C(Y_j)\left[ \Gamma^{-1}_{D}- \widehat{\Gamma}^{-1}_{D}\right] C(Y_j)\Vert_{hs}\\
	&=\Vert\dfrac{1}{n}\sum_{j=1}^{n}C(Y_j)\left[  \widehat{\Gamma}^{-1}_{D}( \widehat{\Gamma}_{D}-\Gamma_D)\Gamma^{-1}_{D}\right] C(Y_j)\Vert_{hs}\\
	&\leq\Vert\dfrac{1}{n}\sum_{j=1}^{n}C(Y_j)\left[  \widehat{\Gamma}^{-1}_{D}( \widehat{\Gamma}_{D}-\Gamma_D)\Gamma^{-1}\right] C(Y_j)\Vert_{hs}\\
	&\leq\Vert \widehat{\Gamma}^{-1}_{D}( \widehat{\Gamma}_{D}-\Gamma_D)\Vert_\infty\dfrac{1}{n}\sum_{j=1}^{n}\Vert C(Y_j)\Vert_{hs}\,\Vert\Gamma^{-1}C(Y_j)\Vert_{hs}.
\end{align*}
Since $\Vert \widehat{\Gamma}_D-\Gamma_D\Vert_\infty=o_p(t_D)$ and $ \Vert \widehat{\Gamma}_{D}^{-1}\Vert_\infty=O_p(\dfrac{1}{t_D})$, we deduce from the preceding inequality that
	$\Vert A_{2n}\Vert_{hs}=o_p(\dfrac{1}{t_D\sqrt{n}})$.
\hfill $\Box$

	\begin{lemma}\label{lem5}
		Under assumptions $(\mathscr{A}_3)$ and $(\mathscr{A}_{10})$,  if we suppose $\Vert \widehat{\Gamma}_{D}-\Gamma_{D}\Vert_\infty=o_p(t_{D}),$ we have
		$$ A_{3n}=\dfrac{1}{n}\sum_{j=1}^{n}C(Y_j) \widehat{\Gamma}^{-1}_{D}C(Y_j)-\dfrac{1}{n}\sum_{j=1}^{n}C_{e_n}(Y_j) \widehat{\Gamma}^{-1}_{D}C_{e_n}(Y_j)=O_p\left( \dfrac{1}{t_{D}\sqrt{n}}\right).$$
	\end{lemma}
\noindent\textit{Proof.}
We can write
	
	\begin{align*}
	A_{3n}&= \dfrac{1}{n}\sum_{j=1}^{n}\left[ C(Y_j)-C_{e_n}(Y_j)\right]  \widehat{\Gamma}^{-1}_{D}C(Y_j)-\dfrac{1}{n}\sum_{j=1}^{n}C_{e_n}(Y_j) \widehat{\Gamma}^{-1}_{D}\left[ C_{e_n}(Y_j)-C(Y_j)\right]\\
	&= A_{31n}-A_{32n}.
	\end{align*}
First, we deal with $A_{31n}$. We have
	\begin{align*}
	A_{31n}&=\dfrac{1}{n}\sum_{j=1}^{n}\left[ R(Y_j)-R_{e_n}(Y_j)\right]  \widehat{\Gamma}^{-1}_{D}C(Y_j)\\
	&+\dfrac{1}{n}\sum_{j=1}^{n}\left[ r_{e_n}(Y_j)\otimes r_{e_n}(Y_j)-r(Y_j)\otimes r(Y_j)\right]  \widehat{\Gamma}^{-1}_{D}C(Y_j)\\
	&=A_{311n}+A_{312n}.
	\end{align*}
Further,
	\begin{align*}
	\sqrt{n}\Vert A_{311n}\Vert_{hs}&=\sqrt{n}\Vert\dfrac{1}{n}\sum_{j=1}^{n}\left[ R(Y_j)-R_{e_n}(Y_j)\right]  \widehat{\Gamma}^{-1}_{D}C(Y_j)\Vert_{hs}\\
	&\leq \dfrac{\Vert \widehat{\Gamma}^{-1}_{D}\Vert_\infty\sqrt{n}}{n}\sum_{j=1}^{n}\Vert M(Y_j)\Vert_{hs}\, \Vert C(Y_j) \Vert_{hs}\,\left|  \dfrac{1}{f(Y_j)}-\dfrac{1}{f_{e_n}(Y_j)}\right|
	\end{align*} 
and since
	\begin{align*}
	\left|  \dfrac{1}{f(Y_j)}-\dfrac{1}{f_{e_n}(Y_j)}\right|&\leq \dfrac{1}{f(Y_j)}\textbf{1}_{\{f(Y_j)<e_n\}}
	\end{align*}
it follows
	\begin{align*}
	\sqrt{n}\Vert A_{311n}\Vert_{hs}&\leq \dfrac{\Vert \widehat{\Gamma}^{-1}_{D}\Vert_\infty\sqrt{n}}{n}\sum_{j=1}^{n}\Vert M(Y_j) \Vert_{hs}\, \Vert C(Y_j)\Vert_{hs}\dfrac{1}{f(Y_j)}\textbf{1}_{\{f(Y_j)<e_n\}} \\
	&= \dfrac{\Vert \widehat{\Gamma}^{-1}_{D}\Vert_\infty\sqrt{n}}{n}\sum_{j=1}^{n}\Vert R(Y_j)\Vert_{hs}\, \Vert C(Y_j)\Vert_{hs}\textbf{1}_{\{f(Y_j)<e_n\}} \\
	&\leq \dfrac{\Vert \widehat{\Gamma}^{-1}_{D}\Vert_\infty\sqrt{n}}{n}\sum_{j=1}^{n}\Vert R(Y_j)\Vert_{hs}\, \Vert R(Y_j)-r(Y_j)\otimes r(Y_j) \Vert_{hs}\textbf{1}_{\{f(Y_j)<e_n\}} \\
	&\leq \dfrac{\Vert \widehat{\Gamma}^{-1}_{D}\Vert_\infty\sqrt{n}}{n}\sum_{j=1}^{n}\left[\Vert R(Y_j)\Vert_{hs}^{2}+  \Vert R(Y_j)\Vert_{hs}\, \Vert r(Y_j)\Vert_{hs}^{2} \right] \textbf{1}_{\{f(Y_j)<e_n\}}.
\end{align*}
Thus
\begin{align*}
	\mathbb{E}\left[\dfrac{\sqrt{n}\Vert A_{311n}\Vert_{hs}}{\Vert \widehat{\Gamma}^{-1}_{D}\Vert_\infty} \right]&\leq \sqrt{n}\mathbb{E}\left[ \Vert R(Y)\Vert_{hs}^{2}\textbf{1}_{\{f(Y)<e_n\}}\right]+\sqrt{n}\mathbb{E}\left[ \Vert R(Y)\Vert_{hs}\, \Vert r(Y)\Vert_{hs}^{2}\textbf{1}_{\{f(Y)<e_n\}}\right] 
\end{align*}
and  since $\Vert \widehat{\Gamma}_{D}^{-1}-\Gamma_{D}^{-1}\Vert_\infty=o_p(t_D)$, $\Vert \widehat{\Gamma}^{-1}_{D}\Vert_\infty=O_p(\dfrac{1}{t_{D}})$,  we deduce from the preceding inequality,  assumption $ (\mathscr{A}_{10})$  and Markov inequality  that  $ A_{311n}=o_p(\dfrac{1}{t_{D}\sqrt{n}})$. On the other hand,
	\begin{align*}
	\sqrt{n}\Vert A_{312n}\Vert_{hs}&\leq \dfrac{\Vert \widehat{\Gamma}_{D}^{-1}\Vert_\infty\sqrt{n}}{n}\sum_{j=1}^{n}\Vert m(Y_j)\otimes m(Y_j)\left[ \dfrac{1}{f^2(Y_j)}-\dfrac{1}{f_{e_n}^{2}(Y_j)}\right]\Vert_{hs}\,\Vert C(Y_j)\Vert_{hs}\\
	&\leq \dfrac{\Vert \widehat{\Gamma}_{D}^{-1}\Vert_\infty\sqrt{n}}{n}\sum_{j=1}^{n}\Vert m(Y_j)\otimes m(Y_j)\Vert_{hs}\,\Vert C(Y_j)\Vert_{hs}\left| \dfrac{1}{f^2(Y_j)}-\dfrac{1}{f_{e_n}^{2}(Y_j)} \right| \\
	&\leq \dfrac{\Vert \widehat{\Gamma}_{D}^{-1}\Vert_\infty\sqrt{n}}{n}\sum_{j=1}^{n}\Vert m(Y_j)\otimes m(Y_j)\Vert_{hs}\,\Vert C(Y_j)\Vert_{hs}\dfrac{1}{f^2(Y_j)}\textbf{1}_{\{f(Y_j)<e_n\}}\\
	&\leq \dfrac{\Vert \widehat{\Gamma}_{D}^{-1}\Vert_\infty\sqrt{n}}{n}\sum_{j=1}^{n}\Vert r(Y_j)\otimes r(Y_j)\Vert_{hs}\,\Vert C(Y_j)\Vert_{hs}\textbf{1}_{\{f(Y_j)<e_n\}}\\
	&\leq \dfrac{\Vert \widehat{\Gamma}_{D}^{-1}\Vert_\infty\sqrt{n}}{n}\sum_{j=1}^{n}\left[ \Vert r(Y_j)\Vert_{H}^{2}\,\Vert R(Y_j)\Vert_{hs}+\Vert r(Y_j)\Vert_{H}^{4}\right] \textbf{1}_{\{f(Y_j)<e_n\}}
\end{align*}
Thus
\begin{align*}
	\mathbb{E}\left[\dfrac{\sqrt{n}\Vert A_{312n}\Vert_{hs}}{\Vert \widehat{\Gamma}_{D}^{-1}\Vert_\infty}\right]&\leq \sqrt{n}\mathbb{E}\left[ \Vert r(Y)\Vert_{H}^{2}\,\Vert R(Y)\Vert_{hs} \textbf{1}_{\{f(Y)<e_n\}} \right]+ \sqrt{n}\mathbb{E}\left[ \Vert r(Y)\Vert_{H}^{4} \textbf{1}_{\{f(Y)<e_n\}} \right]
\end{align*}
	and since $ \Vert \widehat{\Gamma}_{D}^{-1}\Vert_\infty=O_p(\dfrac{1}{t_{D}})$, we deduce from the preceding inequality, assumption $(\mathscr{A}_{10})$ and  Markov inequality  that $ A_{312n}=o_p(\dfrac{1}{t_{D}\sqrt{n}})$. This permits to conclude that  $A_{31n}=A_{311n}+A_{312n}=o_p(\dfrac{1}{t_{D}\sqrt{n}})$. Now, we deal with $A_{32n}$; we have
	\begin{align*}
	A_{32n}&=\dfrac{1}{n}\sum_{j=1}^{n}C_{e_n}(Y_j) \widehat{\Gamma}_{D}^{-1}\left[ R_{e_n}(Y_j)-R(Y_j)\right]-\dfrac{1}{n}\sum_{j=1}^{n}C_{e_n}(Y_j) \widehat{\Gamma}_{D}^{-1}\left[ r_{e_n}(Y_j)\otimes r_{e_n}(Y_j)-r(Y_j)\otimes r(Y_j)\right]\\
	&=A_{321n}-A_{322n}
	\end{align*}
	and
	\begin{align*}
	A_{321n}&=\dfrac{1}{n}\sum_{j=1}^{n}C_{e_n}(Y_j) \widehat{\Gamma}_{D}^{-1}\left[ R_{e_n}(Y_j)-R(Y_j)\right]\\
	&=\dfrac{1}{n}\sum_{j=1}^{n}R_{e_n}(Y_j) \widehat{\Gamma}_{D}^{-1}\left[ R_{e_n}(Y_j)-R(Y_j)\right]-\dfrac{1}{n}\sum_{j=1}^{n}r_{e_n}(Y_j)\otimes r_{e_n}(Y_j) \widehat{\Gamma}_{D}^{-1}\left[ R_{e_n}(Y_j)-R(Y_j)\right]\\
	&=A_{3211n}-A_{3212n}.
	\end{align*}
	Moreover
	\begin{align*}
	\sqrt{n}\Vert A_{3211n}\Vert_{hs}&\leq\dfrac{\Vert \widehat{\Gamma}_{D}^{-1}\Vert_\infty\sqrt{n}}{n}\sum_{j=1}^{n}\Vert R_{e_n}(Y_j)\Vert_{hs}\,\Vert\left[ R_{e_n}(Y_j)-R(Y_j)\right]\Vert_{hs}\\
	&\leq \dfrac{\Vert \widehat{\Gamma}_{D}^{-1}\Vert_\infty\sqrt{n}}{n}\sum_{j=1}^{n}\Vert R_{e_n}(Y_j)\Vert_{hs}\,\Vert m(Y_j)\Vert_{hs}\left| \dfrac{1}{f_{e_n}(Y_j)}-\dfrac{1}{f(Y_j)}\right|  \\
	&\leq \dfrac{\Vert \widehat{\Gamma}_{D}^{-1}\Vert_\infty\sqrt{n}}{n}\sum_{j=1}^{n}\Vert R_{e_n}(Y_j)\Vert_{hs}\,\Vert R(Y_j)\Vert_{hs}\textbf{1}_{\{f(Y_j)<e_n\}}  \\
	&\leq \dfrac{\Vert \widehat{\Gamma}_{D}^{-1}\Vert_\infty\sqrt{n}}{n}\sum_{j=1}^{n}\Vert R(Y_j)\Vert_{hs}^{2}\textbf{1}_{\{f(Y_j)<e_n\}}.
\end{align*}
Hence
\begin{align*}
	\mathbb{E}\left[\dfrac{\sqrt{n}\Vert A_{3211n}\Vert_{hs}}{\Vert \widehat{\Gamma}_{D}^{-1}\Vert_\infty} \right] &\leq \sqrt{n}\mathbb{E}\left[ \Vert R(Y)\Vert_{hs}^{2}\textbf{1}_{\{f(Y)<e_n\}}\right]
\end{align*}
the using    $\Vert \widehat{\Gamma}_{D}^{-1}\Vert_\infty=O_p(\dfrac{1}{t_{D}})$,  assumption $(\mathscr{A}_{10}) $ and Markov inequality, we conclude that $
	A_{3211n}=o_p(\dfrac{1}{t_{D}\sqrt{n}})$. Furthermore, we have
	\begin{align*}
	\sqrt{n}\Vert A_{3212n}\Vert_{hs}&\leq\dfrac{\Vert \widehat{\Gamma}_{D}^{-1}\Vert_\infty\sqrt{n}}{n}\sum_{j=1}^{n}\Vert r(Y_j)\otimes
	r(Y_j)\Vert_{hs}\,\Vert M(Y_j)\Vert_{hs}\left| \dfrac{1}{f_{e_n}(Y_j)}-\dfrac{1}{f(Y_j)}\right| \\
	&\leq\dfrac{\Vert \widehat{\Gamma}_{D}^{-1}\Vert_\infty\sqrt{n}}{n}\sum_{j=1}^{n}\Vert r(Y_j)\Vert_{H}^{2}\,\Vert R(Y_j)\Vert_{hs}\textbf{1}_{\{f(Y_j)<e_n\}}
\end{align*}
what implies
\begin{align*}
	\mathbb{E}\left[ \dfrac{\sqrt{n}\Vert A_{3212n}\Vert_{hs}}{\Vert \widehat{\Gamma}_{D}^{-1}\Vert_\infty}\right] &\leq \sqrt{n}\mathbb{E}\left[\Vert r(Y)\Vert_{H}^{2}\,\Vert R(Y)\Vert_{hs}\textbf{1}_{\{f(Y)<e_n\}} \right]=o_p(1) .
\end{align*}
Thus, $A_{3212n}=o_p(\dfrac{1}{t_{D}\sqrt{n}})$ and 
	we can then conclude that  $A_{321n}=o_p(\dfrac{1}{t_{D}\sqrt{n}}) $. It remains to treat $A_{322n}$. We have:
	\begin{align*}
	A_{322n}&=\dfrac{1}{n}\sum_{j=1}^{n}R_{e_n}(Y_j) \widehat{\Gamma}_{D}^{-1}\left[ r_{e_n}(Y_j)\otimes r_{e_n}(Y_j)-r(Y_j)\otimes r(Y_j)\right]\\
	&-\dfrac{1}{n}\sum_{j=1}^{n}r_{e_n}(Y_j)\otimes r_{e_n}(Y_j) \widehat{\Gamma}_{D}^{-1}\left[ r_{e_n}(Y_j)\otimes r_{e_n}(Y_j)-r(Y_j)\otimes r(Y_j)\right]\\
	&=A_{3221n}-A_{3222n}
	\end{align*}
and since
	\begin{align*}
	\sqrt{n}\Vert A_{3221n}\Vert_{hs}&\leq \Vert \widehat{\Gamma}_{D}^{-1}\Vert_\infty\sqrt{n}\mathbb{E}\left[\Vert r(Y)\Vert_{H}^{2}\,\Vert R(Y)\Vert_{hs}\textbf{1}_{\{f(Y)<e_n\}} \right] 
\end{align*}
we obtain from assumption $(\mathscr{A}_{10})$ that $A_{3221n}=o_p(\dfrac{1}{t_{D}\sqrt{n}})$. Further,	
	\begin{align*}
	\sqrt{n}\mathbb{E}\left[\Vert A_{3222n}\right] \Vert_{hs}&\leq \Vert \widehat{\Gamma}_{D}^{-1}\Vert_\infty \sqrt{n}\mathbb{E}\left[\Vert r(Y)\Vert_{H}^{4}\textbf{1}_{\{f(Y)<e_n\}} \right] 
\end{align*}
and from assumtion $(A_{8})$ and Markov inequality we deduce that $A_{3222n}=o_p(\dfrac{1}{t_{D}\sqrt{n}})$. Consequently, $A_{322n}=o_p(\dfrac{1}{t_{D}\sqrt{n}})$ and  $A_{32n}=o_p(\dfrac{1}{t_{D}\sqrt{n}})$. All of the above permit to conclude that   $A_{3n}=o_p(\dfrac{1}{t_{D}\sqrt{n}})$.
\hfill $\Box$

	\begin{lemma}\label{lem5}
Under assumptions $(\mathscr{A}_3)$, $(\mathscr{A}_6)$ to $(\mathscr{A}_9)$ if we suppose that  $\Vert \widehat{\Gamma}_{D}-\Gamma_{D}\Vert_\infty=o_p(t_{D}),$, we have :
		\begin{align*} 
		A_{4n}&= \dfrac{1}{n}\sum_{j=1}^{n}C_{e_n}(Y_j) \widehat{\Gamma}^{-1}_{D}C_{e_n}(Y_j)-\dfrac{1}{n}\sum_{j=1}^{n} \widehat{C}_{e_n}(Y_j) \widehat{\Gamma}^{-1}_{D} \widehat{C}_{e_n}(Y_j)\\
		&=O_p\left( \dfrac{1}{t_{D}e_n}\left( h^k+\dfrac{\sqrt{ \log(n)}}{h\sqrt{n}}\right) \right) +O_p\left( \dfrac{1}{t_{D}\sqrt{n}}\right)\\
		&=O_p\left( \dfrac{1}{t_{D}n^{\gamma}}\right),
		\end{align*}
where $\gamma$ is a real constant satisfying $0<\gamma<1/4$.
	\end{lemma}
\noindent\textit{Proof.}
	\begin{align*}
	A_{4n}&=\dfrac{1}{n}\sum_{j=1}^{n}\left[ C_{e_n}(Y_j)- \widehat{C}_{e_n}(Y_j)\right]  \widehat{\Gamma}_{D}^{-1}C_{e_n}(Y_j)-\dfrac{1}{n}\sum_{j=1}^{n}\left[  \widehat{C}_{e_n}(Y_j)-C_{e_n}(Y_j)\right] \widehat{\Gamma}_{D}^{-1}\left[  \widehat{C}_{e_n}(Y_j)-C_{e_n}(Y_j)\right]\\
	&+\dfrac{1}{n}\sum_{j=1}^{n}C_{e_n}(Y_j) \widehat{\Gamma}_{D}^{-1}\left[ C_{e_n}(Y_j)- \widehat{C}_{e_n}(Y_j)\right]=A_{41n}-A_{42n}+A_{43n}.
	\end{align*}
	First, we deal with  $A_{41n}$; we have:
		\begin{align*}
		A_{41n}&=\dfrac{1}{n}\sum_{j=1}^{n}\left[ R_{e_n}(Y_j)- \widehat{R}_{e_n}(Y_j)\right] \widehat{\Gamma}_{D}^{-1}C_{e_n}(Y_j)\\
		&+\left[  \widehat{r}_{e_n}(Y_j)\otimes \widehat{r}_{e_n}(Y_j)-r_{e_n}(Y_j)\otimes r_{e_n}(Y_j)\right] \widehat{\Gamma}_{D}^{-1}C_{e_n}(Y_j)\\
		&=\dfrac{1}{n}\sum_{j=1}^{n}\left[ R_{e_n}(Y_j)- \widehat{R}_{e_n}(Y_j)\right] \widehat{\Gamma}_{D}^{-1}C_{e_n}(Y_j)\\
		&+ \dfrac{1}{n}\sum_{j=1}^{n}\left(  \widehat{r}_{e_n}(Y_j)-r_{e_n}(Y_j)\right) \otimes\left(  \widehat{r}_{e_n}(Y_j)-r_{e_n}(Y_j)\right) \widehat{\Gamma}_{D}^{-1}C_{e_n}(Y_j)\\
		&+\dfrac{1}{n}\sum_{j=1}^{n}\left(  \widehat{r}_{e_n}(Y_j)-r_{e_n}(Y_j)\right)\otimes r_{e_n}(Y_j) \widehat{\Gamma}_{D}^{-1}C_{e_n}(Y_j)\\
		&+\dfrac{1}{n}\sum_{j=1}^{n}r_{e_n}(Y_j)\otimes\left(  \widehat{r}_{e_n}(Y_j)-r_{e_n}(Y_j)\right)  \widehat{\Gamma}_{D}^{-1}C_{e_n}(Y_j)\\
		&=A_{411n}+A_{412n}+A_{413n}+A_{414n}
		\end{align*}
		and
		\begin{align*}
		\Vert A_{411n}\Vert_{hs}&\leq\Vert \widehat{\Gamma}_{D}^{-1}\Vert_\infty\dfrac{1}{n}\sum_{j=1}^{n}\Vert R_{e_n}(Y_j)- \widehat{R}_{e_n}(Y_j)\Vert_{hs}\,\Vert C(Y_j)\Vert_{hs}\\
		&\leq \Vert \widehat{\Gamma}_{D}^{-1}\Vert_\infty\dfrac{1}{n}\sum_{j=1}^{n}\Vert R_{e_n}(Y_j)- \widehat{R}_{e_n}(Y_j)\Vert_{hs}\,\Vert C(Y_j)\Vert_{hs}\\
		&\leq \Vert \widehat{\Gamma}_{D}^{-1}\Vert_\infty\dfrac{1}{n}\sum_{j=1}^{n}\Vert\dfrac{ R_{e_n}(Y_j)}{ \widehat{f}_{e_n}(Y_j)}\left[ f(Y_j)- \widehat{f}(Y_j)\right]+\dfrac{1}{ \widehat{f}_{e_n}(Y_j)}\left[  \widehat{M}(Y_j)-M(Y_j)\right] \Vert_{hs}\,\Vert C(Y_j)\Vert_{hs}\\
		&\leq \dfrac{ \Vert \widehat{\Gamma}_{D}^{-1}\Vert_\infty}{e_n}\dfrac{1}{n}\sum_{j=1}^{n}\Vert R(Y_j)\Vert_{hs}\,\Vert C(Y_j)\Vert_{hs}\sup_{y\in \mathbb{R}}| f(y)- \widehat{f}(y)| \\
		&+\dfrac{ \Vert \widehat{\Gamma}_{D}^{-1}\Vert_\infty}{e_n}\dfrac{1}{n}\sum_{j=1}^{n}\,\Vert C(Y_j)\Vert_{hs}\sup_{y\in \mathbb{R}}\Vert M(y)- \widehat{M}(y)\Vert_{hs}.
\end{align*}
It is known from \cite{prakasa}  that
\begin{align}\label{rao}
	\sup_{y\in \mathbb{R}}| \widehat{f}(y)-f(y)|=O_p(h^k+\dfrac{\sqrt{ \log(n)}}{h\sqrt{n}});
	\end{align} 
then, this property together with Lemma \ref{lem2}, assumption $(\mathscr{A}_9)$ and the preceding inequality imply
\begin{align*}
		A_{411n}&=O_p\left( \dfrac{1}{t_{D}e_n}\left( h^k+\dfrac{\sqrt{ \log(n)}}{h\sqrt{n}}\right) \right)=O_p\left( \dfrac{1}{t_{D}n^{\gamma}}\right).
\end{align*}
A similar reasoning, but by using instead of Lemma \ref{lem2} the following result from \cite{yao}:
\begin{align*}
	\sup_{y\in \mathbb{R}}\Vert \widehat{m}(y)-m(y)\Vert_H=O_p\left( h^k+\dfrac{\sqrt{ \log(n)}}{h\sqrt{n}}\right)
	\end{align*}
permits to obtain
		\begin{align*}
		A_{413n}&=O_p\left( \dfrac{1}{t_{D}e_n}\left( h^k+\dfrac{\sqrt{ \log(n)}}{h\sqrt{n}}\right) \right)=O_p\left( \dfrac{1}{t_{D}n^{\gamma}}\right)
\end{align*}
and
\begin{align*}
		A_{414n}&=O_p\left( \dfrac{1}{t_{D}e_n}\left( h^k+\dfrac{\sqrt{ \log(n)}}{h\sqrt{n}}\right) \right)=O_p\left( \dfrac{1}{t_{D}n^{\gamma}}\right).
		\end{align*}
		On the other hand,
		\begin{align*}
		\Vert A_{412n}\Vert_{hs}&=\dfrac{1}{n}\sum_{j=1}^{n}\Vert\left(  \widehat{r}_{e_n}(Y_j)-r_{e_n}(Y_j)\right) \otimes\left(  \widehat{r}_{e_n}(Y_j)-r_{e_n}(Y_j)\right) \widehat{\Gamma}_{D}^{-1}C_{e_n}(Y_j)\Vert_{hs}\\
		&\leq \Vert \widehat{\Gamma}_{D}^{-1}\Vert_\infty\dfrac{1}{n}\sum_{j=1}^{n}\Vert  \widehat{r}_{e_n}(Y_j)-r_{e_n}(Y_j)\Vert_{hs}^{2}\Vert C(Y_j)\Vert_{hs}.
\end{align*}
Similar developpements  as previously done for $\Vert A_{411n}\Vert_{hs}$ permit to obtain $\dfrac{1}{n}\sum_{j=1}^{n}\Vert  \widehat{r}_{e_n}(Y_j)-r_{e_n}(Y_j)\Vert_{hs}^{2}\Vert C(Y_j)\Vert_{hs}=O_p(\dfrac{1}{\sqrt{n}})$, and since $ \Vert \widehat{\Gamma}_{D}^{-1}\Vert_\infty=O_p(\dfrac{1}{t_{D}})$, we conclude that $A_{412n}=O_p\left( \dfrac{1}{t_{D}\sqrt{n}}\right) $. Therefore,
\[	
		A_{41n}=O_p\left( \dfrac{1}{t_{D}e_n}\left( h^k+\dfrac{\sqrt{ \log(n)}}{h\sqrt{n}}\right) \right) +O_p\left( \dfrac{1}{t_{D}\sqrt{n}}\right)=O_p\left( \dfrac{1}{t_{D}n^{\gamma}}\right).
\]
	Further,
	\begin{align*}
	A_{42n}&=\dfrac{1}{n}\sum_{j=1}^{n}\left[  \widehat{R}_{e_n}(Y_j)-R_{e_n}(Y_j)\right] \widehat{\Gamma}_{D}^{-1}\left[  \widehat{R}_{e_n}(Y_j)-R_{e_n}(Y_j)\right]\\
	&-\dfrac{1}{n}\sum_{j=1}^{n}\left[  \widehat{R}_{e_n}(Y_j)-R_{e_n}(Y_j)\right] \widehat{\Gamma}_{D}^{-1}\left[  \widehat{r}_{e_n}(Y_j)\otimes \widehat{r}_{e_n}(Y_j)-r_{e_n}(Y_j)\otimes r_{e_n}(Y_j)\right]\\
	&-\dfrac{1}{n}\sum_{j=1}^{n}\left[  \widehat{r}_{e_n}(Y_j)\otimes \widehat{r}_{e_n}(Y_j)-r_{e_n}(Y_j)\otimes r_{e_n}(Y_j)\right] \widehat{\Gamma}_{D}^{-1}\left[  \widehat{R}_{e_n}(Y_j)-R_{e_n}(Y_j)\right]\\
	&+\dfrac{1}{n}\sum_{j=1}^{n}\left[  \widehat{r}_{e_n}(Y_j)\otimes \widehat{r}_{e_n}(Y_j)-r_{e_n}(Y_j)\otimes r_{e_n}(Y_j)\right] \widehat{\Gamma}_{D}^{-1}\left[  \widehat{r}_{e_n}(Y_j)\otimes \widehat{r}_{e_n}(Y_j)-r_{e_n}(Y_j)\otimes r_{e_n}(Y_j)\right]\\
	&=A_{421n}-A_{422n}-A_{423n}+A_{424n}
	\end{align*}
	and
	\begin{align*}
	\sqrt{n}\Vert A_{421n}\Vert_{hs}&\leq \dfrac{\Vert \widehat{\Gamma}_{D}^{-1}\Vert_\infty\sqrt{n}}{n}\sum_{j=1}^{n}\Vert \widehat{R}_{e_n}(Y_j)- R_{e_n}(Y_j)\Vert_{hs}^{2}\\
	&\leq \dfrac{\Vert \widehat{\Gamma}_{D}^{-1}\Vert_\infty\sqrt{n}}{n}\sum_{j=1}^{n}\Vert\dfrac{ R_{e_n}(Y_j)}{ \widehat{f}_{e_n}(Y_j)}\left[ f_{e_n}(Y_j)- \widehat{f}_{e_n}(Y_j)\right]-\dfrac{1}{ \widehat{f}_{e_n}(Y_j)} \left[  \widehat{M}(Y_j)-M(Y_j)\right]\Vert_{hs}^{2}\\
	&= \dfrac{\Vert \widehat{\Gamma}_{D}^{-1}\Vert_\infty\sqrt{n}}{n}\sum_{j=1}^{n}\Vert\dfrac{ R_{e_n}(Y_j)}{ \widehat{f}_{e_n}(Y_j)}\left[ f_{e_n}(Y_j)- \widehat{f}_{e_n}(Y_j)\right]\Vert_{hs}^{2}\\
	&+\dfrac{\Vert \widehat{\Gamma}_{D}^{-1}\Vert_\infty\sqrt{n}}{n}\sum_{j=1}^{n}\Vert\dfrac{1}{ \widehat{f}_{e_n}(Y_j)} \left[  \widehat{M}(Y_j)-M(Y_j)\right]\Vert_{hs}^{2}\\
	&-\dfrac{2\Vert \widehat{\Gamma}_{D}^{-1}\Vert_\infty\sqrt{n}}{n}\sum_{j=1}^{n}\left\langle \dfrac{ R_{e_n}(Y_j)}{ \widehat{f}_{e_n}(Y_j)}\left[ f_{e_n}(Y_j)- \widehat{f}_{e_n}(Y_j)\right],\dfrac{1}{ \widehat{f}_{e_n}(Y_j)} \left[  \widehat{M}(Y_j)-M(Y_j)\right]\right\rangle _{hs}\\
	&\leq \dfrac{\Vert \widehat{\Gamma}_{D}^{-1}\Vert_\infty\sqrt{n}}{ne_{n}^{2}}\sum_{j=1}^{n}\Vert R_{e_n}(Y_j)\Vert_{hs}^{2}\left| f_{e_n}(Y_j)- \widehat{f}_{e_n}(Y_j)\right|^2\\
	&+\dfrac{\Vert \widehat{\Gamma}_{D}^{-1}\Vert_\infty\sqrt{n}}{ne_{n}^{2}}\sum_{j=1}^{n}\Vert  \widehat{M}(Y_j)-M(Y_j)\Vert_{hs}^{2}\\
	&+\dfrac{2\Vert \widehat{\Gamma}_{D}^{-1}\Vert_\infty\sqrt{n}}{ne_{n}^{2}}\sum_{j=1}^{n}\Vert R_{e_n}(Y_j)\Vert_{hs}\, \Vert  \widehat{M}(Y_j)-M(Y_j)\Vert_{hs}\, \left| f_{e_n}(Y_j)- \widehat{f}_{e_n}(Y_j)\right|\\
	&\leq \dfrac{\Vert \widehat{\Gamma}_{D}^{-1}\Vert_\infty\sqrt{n}}{ne_{n}^{2}}\sum_{j=1}^{n}\Vert R(Y_j)\Vert_{hs}^{2}\left( \sup_{y\in \mathbb{R}}\left|  \widehat{f}(y)-f(y)\right|\right) ^2\\
	&+\dfrac{\Vert \widehat{\Gamma}_{D}^{-1}\Vert_\infty\sqrt{n}}{ne_{n}^{2}}\sum_{j=1}^{n}\left( \sup_{y \in \mathbb{R}}\Vert  \widehat{M}(y)-M(y)\Vert_{hs}\right) ^{2}\\
	&+\dfrac{2\Vert \widehat{\Gamma}_{D}^{-1}\Vert_\infty\sqrt{n}}{ne_{n}^{2}}\sum_{j=1}^{n}\Vert R(Y_j)\Vert_{hs}\, \sup_{y \in \mathbb{R}}\Vert  \widehat{M}(y)-M(y)\Vert_{hs}\, \sup_{y\in \mathbb{R}}\left|  \widehat{f}(y)-f(y)\right|.
	\end{align*}
	From the weak law of large numbers   we obtain $\dfrac{1}{n}\sum_{j=1}^{n}\Vert R(Y_j)\Vert_{hs}=O_p(1)$ and
	$\dfrac{1}{n}\sum_{j=1}^{n}\Vert R(Y_j)\Vert_{hs}^{2}=O_p(1)$. Then using (\ref{rao}), Lemma \ref{lem2}, the assumption   $(\mathscr{A}_{9})$ and the fact that   $\Vert \widehat{\Gamma}_{D}^{-1}\Vert_\infty=O_p(\dfrac{1}{t_{D}})$, we obtain
	\begin{align*}
	\sqrt{n}\Vert A_{421n}\Vert_{hs}&=O_p\left[\dfrac{1}{t_{D}} n^{1/2+2c_2}\left( h^k+\dfrac{1}{h}\sqrt{\dfrac{ \log(n)}{n}}\right)^2\right] \\
	&=O_p\left[ \dfrac{1}{t_{D}}\left( n^{c_2-kc_1+1/4}+n^{c_1+c_2-1/4}\sqrt{ \log(n)}\right)^2\right] \\ 
	&=O_p(\dfrac{1}{t_{D}})
\end{align*}
from what we deduce that $A_{421n}=O_p(\dfrac{1}{t_{D}\sqrt{n}})$. In the same way, we show that $A_{422n}=O_p(\dfrac{1}{t_{D}\sqrt{n}})$, 
	$A_{423n}=O_p(\dfrac{1}{t_{D}\sqrt{n}})$, 
	$A_{424n}=O_p(\dfrac{1}{t_{D}\sqrt{n}})$. Thus, 
	$A_{42n}=O_p(\dfrac{1}{t_{D}\sqrt{n}})$. Now, we deal with 
	$A_{43n}$; since
	  $A_{43n}=(A_{41n})^*$, we  also  have
	\begin{align*}
	A_{43n}&=O_p\left( \dfrac{1}{t_{D}e_n}\left( h^k+\dfrac{\sqrt{ \log(n)}}{h\sqrt{n}}\right) \right) +O_p\left( \dfrac{1}{t_{D}\sqrt{n}}\right)=O_p\left( \dfrac{1}{t_{D}n^{\gamma}}\right)
	\end{align*}
	Finally, we obtain
	\begin{align*} 
	A_{4n}&= \dfrac{1}{n}\sum_{j=1}^{n}C_{e_n}(Y_j) \widehat{\Gamma}_{D}^{-1}C_{e_n}(Y_j)-\dfrac{1}{n}\sum_{j=1}^{n} \widehat{C}_{e_n}(Y_j) \widehat{\Gamma}_{D}^{-1} \widehat{C}_{e_n}(Y_j)\\
	&=O_p\left( \dfrac{1}{t_{D}e_n}\left( h^k+\dfrac{\sqrt{ \log(n)}}{h\sqrt{n}}\right) \right) +O_p\left( \dfrac{1}{t_{D}\sqrt{n}}\right)\\
	&=O_p\left( \dfrac{1}{t_{D}n^{\gamma}}\right).
	\end{align*}
\hfill $\Box$		

\subsection{Proof of the Theorems}
\subsubsection{Proof of Theorem \ref{theo1}}
 Since
\begin{eqnarray*}
& &\mathbb{E} \left[ Var(X|Y)\Gamma^{-1}_{D}Var(X|Y)\right] -\dfrac{1}{n}\sum_{j=1}^{n} \widehat{C}_{e_n}(Y_j) \widehat{\Gamma}^{-1}_{D} \widehat{C}_{e_n}(Y_j)\\
& &\hspace{1cm}=\left( \mathbb{E}\left[ Var(X|Y)\Gamma^{-1}_{D}Var(X|Y)\right]-\dfrac{1}{n}\sum_{j=1}^{n}C(Y_j)\Gamma^{-1}_{D}C(Y_j)\right) \\
& &\hspace{1.5cm}+\left( \dfrac{1}{n}\sum_{j=1}^{n}C(Y_j)\Gamma^{-1}_{D}C(Y_j)-\dfrac{1}{n}\sum_{j=1}^{n}C(Y_j) \widehat{\Gamma}^{-1}_{D}C(Y_j)\right)\\
& &\hspace{1.5cm}+\left( \dfrac{1}{n}\sum_{j=1}^{n}C(Y_j) \widehat{\Gamma}^{-1}_{D}C(Y_j)-\dfrac{1}{n}\sum_{j=1}^{n}C_{e_n}(Y_j) \widehat{\Gamma}^{-1}_{D}C_{e_n}(Y_j)\right) \\
& &\hspace{1.5cm}+\left( \dfrac{1}{n}\sum_{j=1}^{n}C_{e_n}(Y_j) \widehat{\Gamma}^{-1}_{D}C_{e_n}(Y_j)-\dfrac{1}{n}\sum_{j=1}^{n} \widehat{C}_{e_n}(Y_j) \widehat{\Gamma}^{-1}_{D} \widehat{C}_{e_n}(Y_j)\right) \\
& &\hspace{1cm}=A_{1n}+A_{2n}+A_{3n}+A_{4n}.
\end{eqnarray*}
From the central limit theorem we have $A_{1n}=O_p\left(\frac{1}{\sqrt{n}}\right)$; then  the  required result is obtained by applying lemmas \ref{lem3} to \ref{lem5}.
\subsubsection{Proof of Theorem \ref{theo2}}
Putting 
\begin{equation}\label{g}
G=2\Gamma_e+\Psi-\Gamma,\,\,\,  G_D=2\Gamma_e+\mathbb{E}[Var(X|Y)\Gamma^{-1}_{D}Var(X|Y)]-\Gamma
\end{equation} 
and
\begin{equation}\label{gchap}
 \widehat{G}=2 \widehat{\Gamma}_{e,n}+\dfrac{1}{n}\sum\limits_{j=1}^{n} \widehat{C}_{e_n}(Y_j) \widehat{\Gamma}_{D}^{-1} \widehat{C}_{e_n}(Y_j)-\Gamma_n=2 \widehat{\Gamma}_{e,n}+ \widehat{\Psi}_{e_n,D}-\Gamma_n,
\end{equation}
we have:
\begin{eqnarray}\label{ineq1}
\Vert \widehat{\Gamma}_{I,n}-\Gamma_{I}\Vert_{hs}
& \leq & \Vert\Gamma^{-1}G-\Gamma^{-1}G_D\Vert_{hs}+\Vert\Gamma^{-1}G_D-\Gamma^{-1}_{D}G_D\Vert_{hs}+\Vert\Gamma^{-1}_{D}G_D- \widehat\Gamma^{-1}_{D}G_D\Vert_{hs}\nonumber\\
& &+\Vert \widehat\Gamma^{-1}_{D}G_D- \widehat\Gamma^{-1}_{D} \widehat{G}\Vert_{hs}\nonumber\\
&=&K_{1n}+K_{2n}+K_{3n}+K_{4n}.
\end{eqnarray}
First, 
\[
\lim_{D\rightarrow +\infty}K_{1n}=\lim_{D\rightarrow +\infty} 
\Vert\mathbb{E}[\Gamma^{-1}Var(X|Y)(\Gamma^{-1}-\Gamma^{-1}_{D})Var(X|Y)]\Vert_{hs}=0
\]
and since $K_{2n}=\Vert\Gamma^{-1}G_D-\Gamma^{-1}_{D}G_D\Vert_{hs}
\leq \Vert(\Gamma^{-1}-\Gamma^{-1}_{D})G\Vert_{hs}$, we also have $\lim_{D\rightarrow +\infty}K_{2n}=0$. Further,
\begin{eqnarray*}
K_{3n}&=&\Vert\Gamma^{-1}_{D}G_D- \widehat\Gamma^{-1}_{D}G_D\Vert_{hs}
\leq \Vert(\Gamma^{-1}_{D}- \widehat\Gamma^{-1}_{D})G\Vert_{hs}\\
&\leq& \Vert \widehat\Gamma^{-1}_{D}( \widehat\Gamma_D-\Gamma_D)\Gamma^{-1}_{D}G\Vert_{hs}
\leq \Vert \widehat\Gamma^{-1}_{D}\Vert_{hs}\,\Vert \widehat\Gamma_D-\Gamma_D\Vert_{hs}\,\Vert\Gamma^{-1}G\Vert_{hs},
\end{eqnarray*}
then since $\Vert \widehat\Gamma^{-1}_{D}\Vert_{hs}=O_p(1/t_D)$ and $\Vert \widehat\Gamma_D-\Gamma_D\Vert_{hs}=o_p(t_D)$, we deduce that $ K_{3n}=o_p(1)$. On the other hand
\begin{eqnarray*}
K_{4n}&=&\Vert \widehat\Gamma^{-1}_{D}G_D- \widehat\Gamma^{-1}_{D} \widehat{G}\Vert_{hs}\\
&\leq& \Vert \widehat\Gamma^{-1}_{D}\left(  \widehat{\Gamma}_{e,n}-\Gamma_e\right) \Vert_{hs}+\Vert \widehat\Gamma^{-1}_{D}\left\lbrace\mathbb{E}[Var(X|Y)\Gamma^{-1}_{D}Var(X|Y)]-\Psi_{e,n} \right\rbrace \Vert_{hs}\\
& &+\Vert \widehat\Gamma^{-1}_{D}(\Gamma_n-\Gamma)\Vert_{hs}.
\end{eqnarray*}
Since  $\Vert \widehat{\Gamma}_{e,n}-\Gamma_e\Vert_{hs}=O_p(1/\sqrt{n})$ (see \cite{ferre}), we deduce from the preceding inequality that 
\[
K_{4n}=O_p\left( \dfrac{1}{t_D\sqrt{n}}\right) +O_p\left( \dfrac{1}{t_{D}^{2}n^{\gamma}}\right) =o_p(1).
\]
Then using (\ref{ineq1})  and the previous results we obtain:
$ \widehat{\Gamma}_{I,n}-\Gamma_{I}=o_p(1)$.
\subsubsection{Proof of Theorem \ref{theo3}}
Denoting by  $( \widehat{\beta}_{k})_{1\leq k\leq K}$ the orthonormal eigenvectors associated with the $K$ largest eigenvalues $ \widehat{\lambda}_{1}> \widehat{\lambda}_{2}>\cdots> \widehat{\lambda}_{K}>0$ of  $ \widehat\Gamma_{D}^{-1} \widehat{G}$ and by $(\beta_k)_{1\leq k\leq K}$ the orthonormal eigenvectors associated with the $K$ largest eigenvalues $\lambda_{1}>\lambda_{2}>\cdots>\lambda_{K}>0$ of   $\Gamma^{-1}G$, where $G$ and $\widehat{G}$ are defined in (\ref{g}) and (\ref{gchap}), we will only show the convergence of the vector $ \widehat{\beta}_{1}$  as the proof for the others are the same. Clearly, $\beta_1=\lambda_1^{-1}\Gamma_2\eta$ and $\widehat{\beta}_{1}=\widehat{\lambda}_{1}^{-1} \widehat{\Gamma}_2 \widehat{\eta}$ where $\Gamma_2=\Gamma^{-1}\{2\Gamma_e+\Psi-\Gamma\}\Gamma^{-1/2}$, $\eta=\Gamma^{1/2}\beta_1$ and $ \widehat{\eta}= \widehat{\Gamma}_{D}^{1/2} \widehat{\beta}_{1}$ with 
\[	
	 \widehat{\Gamma}_2\dfrac{1}{ \widehat{\lambda}_{1}} \widehat{\Gamma}_{D}^{-1}\left[2 \widehat{\Gamma}_{e,n}+ \widehat{\Psi}_{e_n}-\Gamma_n\right] \widehat{\Gamma}_{D}^{-1/2}.
\]
Hence	
\begin{eqnarray*}
\Vert \widehat{\beta}_{1}-\beta_1\Vert_{H}
	&\leq & \Vert\dfrac{1}{ \widehat{\lambda}_{1}}( \widehat\Gamma_2-\Gamma_2) \widehat\eta\Vert_{H}+\Vert\dfrac{1}{\lambda_1}\Gamma_2( \widehat{\eta}-\eta)\Vert_{H}+\Vert(\dfrac{1}{ \widehat{\lambda}_{1}}-\dfrac{1}{\lambda_1})\Gamma_2 \widehat{\eta}\Vert_{H}\\
	&\leq& \dfrac{\Vert \widehat\eta\Vert_{H}}{| \widehat{\lambda}_{1}|}\Vert \widehat\Gamma_2-\Gamma_2\Vert_\infty+\dfrac{\Vert\Gamma_2\Vert_\infty}{|\lambda_1|}\,\Vert \widehat{\eta}-\eta\Vert_{H}+\dfrac{| \widehat{\lambda}_{1}-\lambda_1|}{|\lambda_1 \widehat{\lambda}_{1}|}\Vert\Gamma_2 \widehat{\eta}\Vert_{H}.
\end{eqnarray*}
Then from Lemma 1 in \cite{ferre} we obtain the inequalities
\begin{equation}\label{lem1yao}
	| \widehat{\lambda}_{1}-\lambda_1|\leq \Vert \widehat{\Gamma}_{D}^{1/2} \widehat\Gamma_2-\Gamma^{1/2}\Gamma_2\Vert_\infty\,\,\,
\textrm{ and }\,\,\,
	\Vert \widehat{\eta}-\eta\Vert_{H}\leq C_9\Vert \widehat{\Gamma}_{D}^{1/2} \widehat\Gamma_2-\Gamma^{1/2}\Gamma_2\Vert_\infty, 
\end{equation}
where $C_9$ is an appropriate positive constant. Then, putting $L_n=\Vert \widehat\Gamma_2-\Gamma_2\Vert_\infty$ and $
M_n=\Vert\Gamma_{n}^{1/2} \widehat\Gamma_2-\Gamma^{1/2}\Gamma_2\Vert_\infty$, we have
\begin{equation}\label{beta}
	\Vert \widehat\beta_{1}-\beta_1\Vert_{H}\leq\dfrac{\Vert \widehat\eta\Vert_{H}}{| \widehat{\lambda}_{1}|}L_n+\left( C_{10}+\dfrac{C_{11}\Vert \widehat\eta\Vert_{H}}{| \widehat{\lambda}_{1}|}\right)\, M_n.
\end{equation}
	Let us verify that $L_n=o_p(1)$ and $M_n=o_p(1)$. First,
\begin{eqnarray*}
	L_n&=&\Vert\Gamma^{-1}G\Gamma^{-1/2}- \widehat{\Gamma}_{D}^{-1} \widehat{G} \widehat{\Gamma}_{D}^{-1/2}\Vert_\infty\\
	&\leq& \Vert\Gamma^{-1}G\Gamma^{-1/2}-\Gamma_{D}^{-1}G\Gamma_{D}^{-1/2}\Vert_\infty+\Vert\Gamma^{-1}_{D}G\Gamma_{D}^{-1/2}- \widehat{\Gamma}_{D}^{-1}G \widehat{\Gamma}_{D}^{-1/2}\Vert_\infty+\Vert \widehat{\Gamma}_{D}^{-1}(G- \widehat{G}) \widehat{\Gamma}_{D}^{-1/2}\Vert_\infty\\
	&=&L_{1n}+L_{2n}+L_{3n}.
	\end{eqnarray*}
	We know that  $\Gamma_{D}^{-1}G\Gamma_{D}^{-1/2}=\Pi_{D}\Gamma^{-1}G\Gamma^{-1/2}\Pi_{D}$ and putting
	$\Pi_{D}^{\bot }=I-\Pi_{D}$ we have  
	\begin{align*}
	\Gamma_{D}^{-1}G\Gamma_{D}^{-1/2}-\Gamma^{-1}G\Gamma^{-1/2}&=\Pi_{D}\Gamma^{-1}G\Gamma^{-1/2}\Pi_{D}-\Gamma^{-1}G\Gamma^{-1/2}\\
	&=\Pi_{D}\Gamma^{-1}G\Gamma^{-1/2}\Pi_{D}-\Pi_{D}^{\bot}\Gamma^{-1}G\Gamma^{-1/2}-\Pi_{D}\Gamma^{-1}G\Gamma^{-1/2}\\
	&=\Pi_{D}\Gamma^{-1}G\Gamma^{-1/2}-\Pi_{D}\Gamma^{-1}G\Gamma^{-1/2}\Pi_{D}^{\bot}-\Pi_{D}^{\bot}\Gamma^{-1}G\Gamma^{-1/2}\\
	&\noindent-\Pi_{D}\Gamma^{-1}G\Gamma^{-1/2}\\
	&=-\Pi_{D}\Gamma^{-1}G\Gamma^{-1/2}\Pi_{D}^{\bot}-\Pi_{D}^{\bot}\Gamma^{-1}G\Gamma^{-1/2}.
\end{align*}
Thus
\begin{align*}
	L_{1n}
	&\leq\Vert\Pi_{D}\Gamma^{-1}G\Gamma^{-1/2}\Pi_{D}^{\bot}\Vert_\infty+\Vert\Pi_{D}^{\bot}\Gamma^{-1}G\Gamma^{-1/2}\Vert_\infty\\
	&\leq\Vert\Gamma^{-1}G\Gamma^{-1/2}\Pi_{D}^{\bot}\Vert_\infty+\Vert\Pi_{D}^{\bot}\Gamma^{-1}G\Gamma^{-1/2}\Vert_\infty
\end{align*}
and, consequently, $\lim_{D\rightarrow +\infty}L_{1n}=0$ because  $\lim_{D\rightarrow +\infty }\Pi_{D}^{\bot}=0$. On the other hand,
	\begin{align*}
	L_{2n}&\leq\Vert(\Gamma^{-1}_{D}- \widehat{\Gamma}_{D}^{-1})G\Gamma^{-1/2}_{D}\Vert_\infty+\Vert\Gamma^{-1}_{D}G(\Gamma^{-1/2}_{D}- \widehat{\Gamma}_{D}^{-1/2})\Vert_\infty+\Vert(\Gamma^{-1}_{D}- \widehat{\Gamma}_{D}^{-1})G(\Gamma^{-1/2}_{D}- \widehat{\Gamma}_{D}^{-1/2})\Vert_\infty
	\end{align*}
	and
	\begin{align*}
	\Vert(\Gamma^{-1}_{D}- \widehat{\Gamma}_{D}^{-1})G\Gamma^{-1/2}_{D}\Vert_\infty&=\Vert \widehat{\Gamma}_{D}^{-1}( \widehat{\Gamma}_{D}-\Gamma_D)\Gamma^{-1}_{D}G\Gamma^{-1/2}_{D}\Vert_\infty\\
	&\leq \Vert \widehat{\Gamma}_{D}^{-1}( \widehat{\Gamma}_{D}-\Gamma_D)\Gamma^{-1}G\Gamma^{-1/2}\Vert_\infty\\
	&\leq \Vert \widehat{\Gamma}_{D}^{-1}\Vert_\infty\,\Vert \widehat{\Gamma}_{D}-\Gamma_D\Vert_\infty\,\Vert\Gamma^{-1}G\Gamma^{-1/2}\Vert_\infty\\
	&=O_p\left( \dfrac{1}{t_D\sqrt{n}}\right) \\
	&=o_p(1).
	\end{align*}
	Using the following properties of operators (see, e.g., \cite{fukumizu}): 
\[
A^{-1/2}-B^{-1/2}=A^{-1/2}(B^{3/2}-A^{3/2})B^{-3/2}+(A-B)B^{-3/2}\,\,\,\textrm{ and }\,\,\,\Vert A^{3/2}-B^{3/2}\Vert_\infty\leq C_{12}\Vert A-B\Vert_\infty
\]
we obtain:
	\[
	\Vert\Gamma^{-1}_{D}G(\Gamma^{-1/2}_{D}- \widehat{\Gamma}_{D}^{-1/2})\Vert_\infty=O_p\left( \dfrac{1}{t_{D}^{3/2}\sqrt{n}}\right)
\]
 and 
\[
\Vert(\Gamma^{-1}_{D}- \widehat{\Gamma}_{D}^{-1})G(\Gamma^{-1/2}_{D}- \widehat{\Gamma}_{D}^{-1/2})\Vert_\infty=O_p\left( \dfrac{1}{t_{D}^{5/2}n}\right).
\]
Therefore, $L_{2n}=o_p(1)$. For dealing with the last term  $L_{3n}$ we consider the operator $G_D=2\Gamma_e+\mathbb{E}[Var(X|Y)\Gamma_{D}^{-1}Var(X|Y)]-\Gamma$  and we have 
	\begin{align*}
	\Vert \widehat{\Gamma}_{D}^{-1}(G-G_D) \widehat{\Gamma}_{D}^{-1/2}\Vert_{hs}&\leq 	\Vert( \widehat{\Gamma}_{D}^{-1}-\Gamma_{D}^{-1})(G-G_D)( \widehat{\Gamma}_{D}^{-1/2}-\Gamma_{D}^{-1/2})\Vert_{hs}\\
&+\Vert( \widehat{\Gamma}_{D}^{-1}-\Gamma_{D}^{-1})(G-G_D)\Gamma_{D}^{-1/2}\Vert_{hs}\\
	&+\Vert\Gamma_{D}^{-1}(G-G_D)( \widehat{\Gamma}_{D}^{-1/2}-\Gamma_{D}^{-1/2})\Vert_{hs}+\Vert\Gamma_{D}^{-1}(G-G_D)\Gamma_{D}^{-1/2}\Vert_{hs}\\
	&\leq 	\Vert \widehat{\Gamma}_{D}^{-1}-\Gamma_{D}^{-1}\Vert_{hs}\,\Vert G-G_D\Vert_{hs}\,\Vert \widehat{\Gamma}_{D}^{-1/2}-\Gamma_{D}^{-1/2}\Vert_{hs}\\
&+\Vert \widehat{\Gamma}_{D}^{-1}-\Gamma_{D}^{-1}\Vert_{hs}\,\Vert(G-G_D)\Gamma^{-1/2}\Vert_{hs}\\
	&+\Vert\Gamma^{-1}(G-G_D)\Vert_{hs}\,\Vert \widehat{\Gamma}_{D}^{-1/2}-\Gamma_{D}^{-1/2}\Vert_{hs}+\Vert\Gamma^{-1}(G-G_D)\Gamma^{-1/2}\Vert_{hs}\\
	 &\leq\Vert \widehat{\Gamma}_{D}^{-1}-\Gamma_{D}^{-1}\Vert_{hs}\,\Vert G\Vert_{hs}\,\Vert \widehat{\Gamma}_{D}^{-1/2}-\Gamma_{D}^{-1/2}\Vert_{hs}+\Vert \widehat{\Gamma}_{D}^{-1}-\Gamma_{D}^{-1}\Vert_{hs}\,\Vert G\Gamma^{-1/2}\Vert_{hs}\\
	 &+\Vert\Gamma^{-1}G\Vert_{hs}\,\Vert \widehat{\Gamma}_{D}^{-1/2}-\Gamma_{D}^{-1/2}\Vert_{hs}+\Vert\Gamma^{-1}(G-G_D)\Gamma^{-1/2}\Vert_{hs}\\
	 &=O_p\left( \dfrac{1}{t_D\sqrt{n}}\right) +O_p\left( \dfrac{1}{t_{D}^{3/2}\sqrt{n}}\right) +O_p\left( \dfrac{1}{t_{D}^{5/2}n}\right)+o_p(1)\\
	 &=o_p(1).
	\end{align*}
Thus
	\begin{align*}
	L_{3n}&\leq \Vert \widehat{\Gamma}_{D}^{-1}(G-G_D) \widehat{\Gamma}_{D}^{-1/2}\Vert_{hs}+\Vert \widehat{\Gamma}_{D}^{-1}(G_D- \widehat{G}) \widehat{\Gamma}_{D}^{-1/2}\Vert_{hs}\\
	&\leq \Vert \widehat{\Gamma}_{D}^{-1}(\Gamma-\Gamma_n) \widehat{\Gamma}_{D}^{-1/2}\Vert_\infty+2\Vert \widehat{\Gamma}_{D}^{-1}\left( \Gamma_e- \widehat{\Gamma}_{e,n}\right)  \widehat{\Gamma}_{D}^{-1/2}\Vert_\infty\\
	&+\Vert \widehat{\Gamma}_{D}^{-1}\left(\mathbb{E}\left[ Var(X|Y)\Gamma^{-1}_{D}Var(X|Y)\right]-\dfrac{1}{n}\sum_{j=1}^{n} \widehat{C}_{e_n}(Y_j) \widehat{\Gamma}_{D}^{-1} \widehat{C}_{e_n}(Y_j)\right)  \widehat{\Gamma}_{D}^{-1/2}\Vert_\infty+o_p(1)\\
	&=O_p\left(\dfrac{1}{t_{D}^{3/2}\sqrt{n}} \right) + O_p\left(\dfrac{1}{t_{D}^{5/2}n^{\gamma}} \right)+o_p(1)\\
	&=o_p(1)
	\end{align*}
	From all what precedes we deduce that $L_{n}=o_p(1)$. From  similar reasoning we also obtain   $M_n=o_p(1)$. Then from (\ref{lem1yao}) and what precedes, we deduce that 
$\widehat{\lambda}_{1}^{-1}=O_p(1)$ and $\Vert \widehat{\eta}\Vert_{H}=O_p(1)$. Therefore,   (\ref{beta}) allows to conclude that $\Vert \widehat\beta_{1}-\beta_1\Vert_{H}=o_p(1)$.

\end{document}